\documentclass[a4paper,11pt]{amsproc}
\usepackage{t1enc,amssymb,epsfig}
\usepackage[english]{babel}

\input xy
\xyoption{all}

\newtheorem{proposition}{Proposition}[section]
\newtheorem{theorem}[proposition]{Theorem}
\newtheorem{lemma}[proposition]{Lemma}
\newtheorem{definition}[proposition]{Definition}

\newcommand{\C}{\mathbb{C}}
\newcommand{\Z}{\mathbb{Z}}

\newcommand{\R}{\mathbb{R}}

\begin{document}

\title[Introduction to the resolution of
    singularities]{\textbf{Introduction to  
        Jung's method of
        resolution of singularities}} 
\author{\sc Patrick Popescu-Pampu}
\address{Universit{\'e} Lille 1, UFR de Maths., B\^atiment M2\\
     Cit\'e Scientifique, 59655, Villeneuve d'Ascq Cedex, France.}
   \email{patrick.popescu@math.univ-lille1.fr}


\medskip

\subjclass[2000]{Primary 32 S 45; Secondary 32 C 20.}

This paper appeared in \textit{Topology of Algebraic Varieties and Singularities.} 
  Proceedings of the conference in honor of Anatoly Libgober's 60th birthday. 
J. I. Cogolludo-Agustin and E. Hironaka eds. Contemporary Mathematics {\bf 538} 
AMS, 2011, 401-432.

\medskip

\maketitle
\thispagestyle{empty}

\begin{quote}

{\small \emph{Socrates}: Now, Glaucon, let's think 
  about the ignorance of human beings and 
  their education in the form of an
   allegory. Imagine them living underground in a 
   kind of cave....
   They see only the     
   shadows the light from the fire throws on  
   the wall of the cave       
   in front of them....  So, it's obvious that 
   for these prisoners 
   the truth would be no more than the shadows of 
   objects.... Now
   let's consider how they might be released and 
   cured of their
   ignorance. Imagine that one man is set free and 
   forced to turn
   around and walk toward the light. Looking at the 
   light will be painful....

   \hfill (Plato: \emph{The Republic.} Book Seven)
}

\end{quote}

\par\medskip\centerline{\rule{2cm}{0.2mm}}\medskip
\setcounter{section}{0}

\tableofcontents

\section{Introduction}

The present notes originated in the 
introductory course given at the  
Trieste Summer School on Resolution of Singularities, 
in June 2006. They 
focus on the resolution of complex   
analytic curves and surfaces by Jung's method. They do
not contain detailed proofs, but mainly explanations 
of the central concepts and of their 
interrelations, as well as heuristics. 
\medskip

If we have begun this text by a famous quotation 
from Plato, it is because we
believe that the citation is related to the general 
philosophical idea of 
\emph{resolution of singularities}. This idea 
corresponds  
to a \emph{desire} which one can search to fulfil 
in various contexts, 
mathematics being one of them. And inside mathematics, 
one can search 
to fulfil it inside various of its branches, 
for various categories of objects.  
A formulation of the desire could be: 
\begin{equation}  \label{fundesire}
    \begin{array}{c}
         \mbox{\emph{Given a complicated object, 
      represent it as the image}} \\  
          \mbox{\emph{of a less   complicated one.} }
     \end{array}   
\end{equation}

As a very simple example, consider a well-known 
problem of
elementary topology. It asks to prove that, 
given three villages
and three wells, one cannot construct roads
joining each village to each well, and such that 
the roads meet only
at the villages or the wells. While dealing with 
the problem, one is 
naturally led to construct diagrams as in Figure 1. 
In it, there is an undesired  
crossing at the point $P$. One has to prove that 
such crossings are
unavoidable. Or, reformulating the problem, 
that one cannot embed  the
abstract graph  drawn in Figure 1 into the plane. 
In this drawing one
has \emph{to imagine that the  crossing point $P$  
is not present.}

  {\tt    \setlength{\unitlength}{0.92pt}}
 \begin{figure}
  \epsfig{file=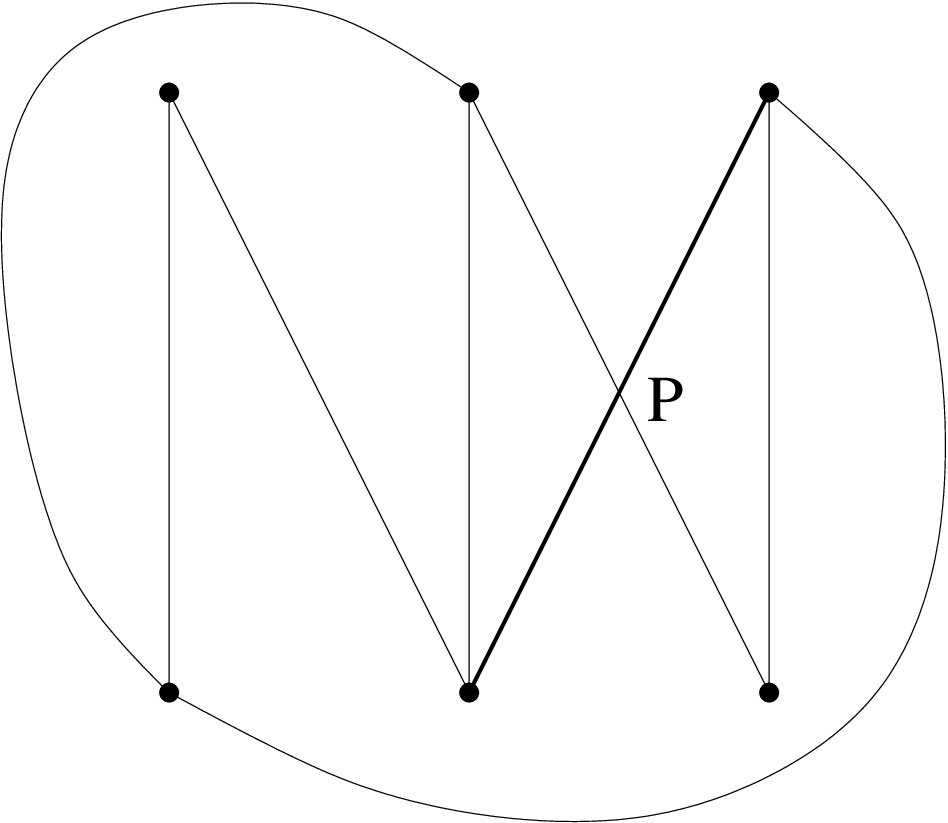, height= 30 mm}
   \caption{An undesired crossing}
   \end{figure}

Now, a little introspection shows that we imagine 
very easily such an \emph{elimination  
of crossing points}, that this is in fact part of 
our automatic
toolkit for understanding images. For example, 
think  of the
perspective drawing of a cube, in which the 
3-dimensional object  
``jumps to the eyes''. 

In the previous examples, we can isolate an elementary 
operation of local 
nature, in which we imagine two lines crossing on 
the plane as the
projections of two uncrossing  
lines in space. This is the easiest example of 
\emph{resolution of  singularities}. The singularity  
which is being resolved is the  germ of the union of 
the two lines at
their meeting point $P$, their  
resolution is the union of the germs of the space 
lines at the preimages of the  point $P$.  

Let us give now a 2-dimensional example. 
We start again with a topological object, the  
\emph{projective plane}. One can present it to 
beginners as the  
surface obtained from a disc by identifying 
opposite points of its boundary.  A way to do   
this identification in 3-dimensional space is 
to divide first the boundary into four equal  
arcs; secondly, to deform the disc till one glues 
two opposite arcs along
a segment;   finally to glue the two remaining 
arcs along the same
segment. One gets like this a  so-called 
\emph{cross-cap} (see Figure 2). 

{\tt    \setlength{\unitlength}{0.92pt}}
 \begin{figure} 
\epsfig{file=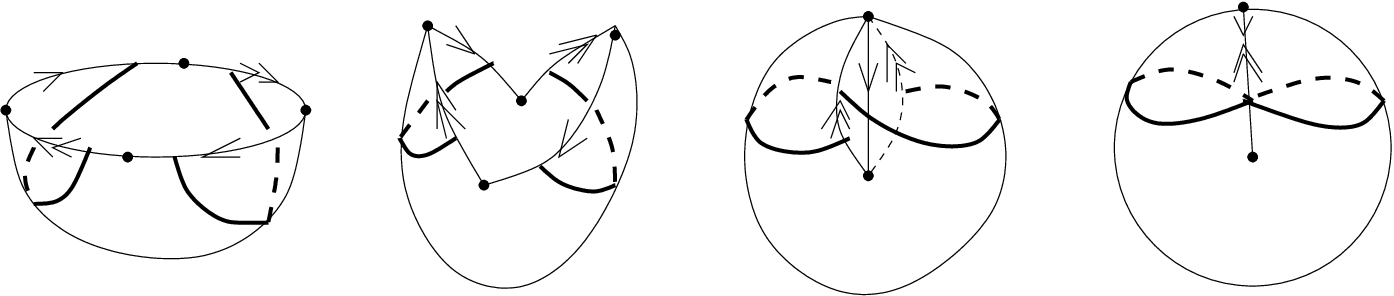, height= 27 mm}
 \caption{Construction of a cross-cap}
 \end{figure}

Arrived at this point, one has to eliminate 
by the imagination this
segment of self-crossing of the  
surface in space, in order to get topologically 
the projective plane. This is an exercise  
analogous to  the one performed before with 
the graph of Figure 1 
in order to get rid of the crossing point $P$. But
now it is more difficult to  
imagine it, as in everyday life we do not interpret 
self-crossing surfaces in space as  
projections of surfaces in a higher-dimensional space. 
In fact, the need to do such  
interpretations was historically one of the 
driving forces of the
elaboration of a mathematical theory of 
higher-dimensional spaces.

For the moment we have constructed the cross-cap 
only as a topological
space. But one can give  
an algebraic model of it. Let us cite Hilbert 
\& Cohn-Vossen
\cite[VI.47]{HCV 90}, who explain a particularly 
nice way to get  a defining equation: 

\begin{quote}

{\small   There is an algebraic surface of this form. 
Its equation is
   \begin{equation} \label{eqcrosscap}
      (k_1x^2 + k_2 y^2)(x^2+y^2+z^2)-2z(x^2+y^2)=0.
   \end{equation}
   This surface is connected with a construction 
   in differential
   geometry. On any surface $F$,  
   we begin with a point $P$ at which the curvature 
   of $F$ is
   positive. Then we construct all the circles  
   of normal curvature at $P$.    
   This family of circles sweeps out [a
   cross-cap], where the line of  
   self-intersection is a segment of the normal 
   to the surface $F$ at
   $P$. The equation given  
   above is referred to the rectangular 
   coordinate system with $P$ as
   origin and with the principal  
   directions of $F$ at the point $P$ as $x$-axis 
   and $y$-axis.  
   $k_1$ and $k_2$ are the principal  
   curvatures of $F$ at the point $P$.}
\end{quote}

{\tt    \setlength{\unitlength}{0.92pt}}
 \begin{figure} 
 \epsfig{file=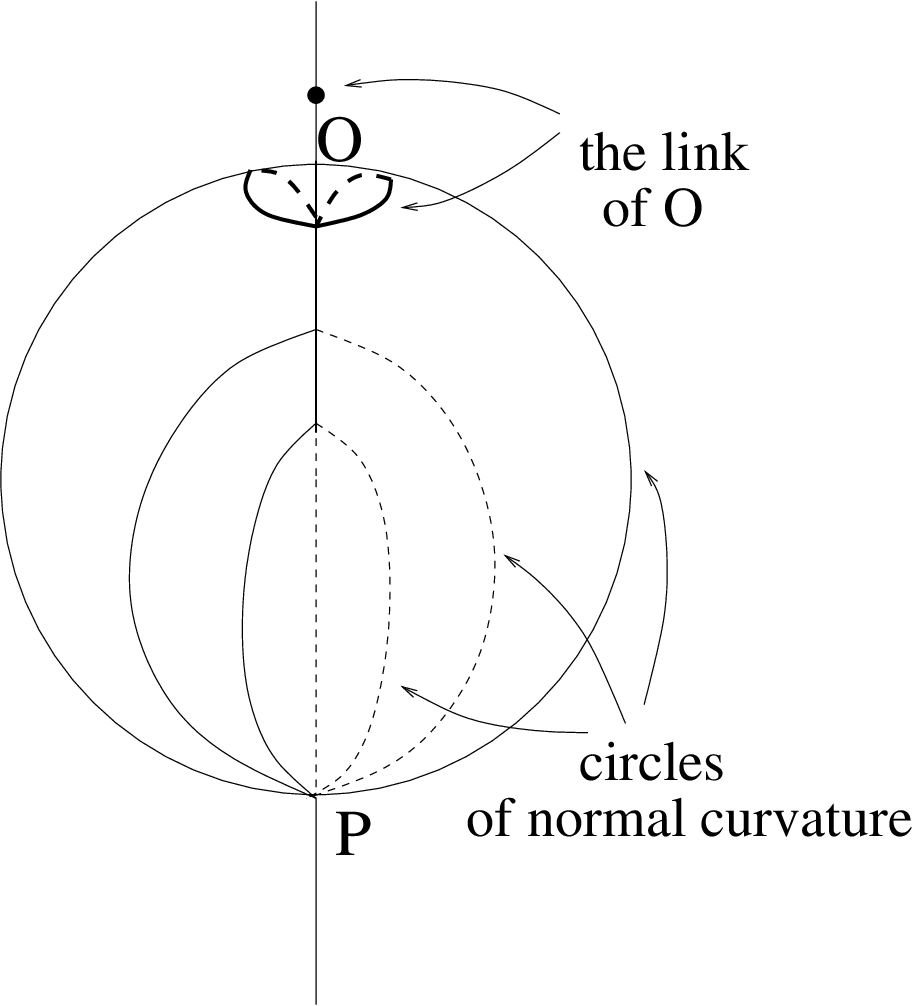, height= 60 mm}
  \caption{Why a cross-cap cannot be algebraic}
 \end{figure}

Something is subtly wrong in the previous statements: 
the locus of
points satisfying    
equation (\ref{eqcrosscap}) is not reduced to 
the cross-cap, but it
contains also the entire  
axis of the variable $z$. Notice that the 
intersection of 
the cross-cap and of the $z$-axis is equal to 
the segment of
self-intersection of the cross-cap (see Figure 3).  
This shows that the cross-cap described by Hilbert 
and Cohn-Vossen is
in fact \emph{a real semi-algebraic 
  surface}, which means that  
it can be defined by a finite set of polynomial 
equations and
\emph{inequations}. At this point, we  
cannot resist the temptation of explaining why no 
(abstract) real algebraic surface
can be homeomorphic to a cross-cap. This is 
a consequence of the
following theorem of Sullivan \cite{S 71}:

\begin{theorem}
     Let $X$ be a real algebraic set and $O$ 
     a point of $X$. Then the
     link of $O$ in $X$ has  
     even Euler characteristic.
 \end{theorem}
 
 The \emph{link} of $O$ in $X$ is, by definition, 
the boundary of a regular
 neighborhood of $O$ in $X$,  obtained by 
 intersecting  $X$ with a sufficiently small ball  
 centered at $O$, after having embedded $X$ in 
an ambient euclidean
 space. In our case, consider as a point $O$ 
the opposite of  $P$ on the cross-cap.  
 The link of $O$ is homeomorphic to 
 an $\infty$-{\em shaped} curve,
 whose Euler  
 characteristic is equal to $-1$. 
 By Sullivan's theorem, this shows
 that a real algebraic surface  
 containing the cross-cap must contain 
 also other points in the
 neighborhood of $O$. In this way,  
 one understands better the presence of 
 the ``stick'' getting out of $O$. 
 
 There is a famous surface in the real 
 3-dimensional space, whose topology 
 captures precisely the local topology of the surface
 (\ref{eqcrosscap}) in the  
 neighborhood of the point $O$. 
  It is called \emph{Whitney's
   umbrella}, and is defined by the equation: 
 \begin{equation} \label{Whitneyeq}
      x^2-zy^2=0.
 \end{equation} 

{\tt    \setlength{\unitlength}{0.92pt}}
 \begin{figure} 
  \epsfig{file=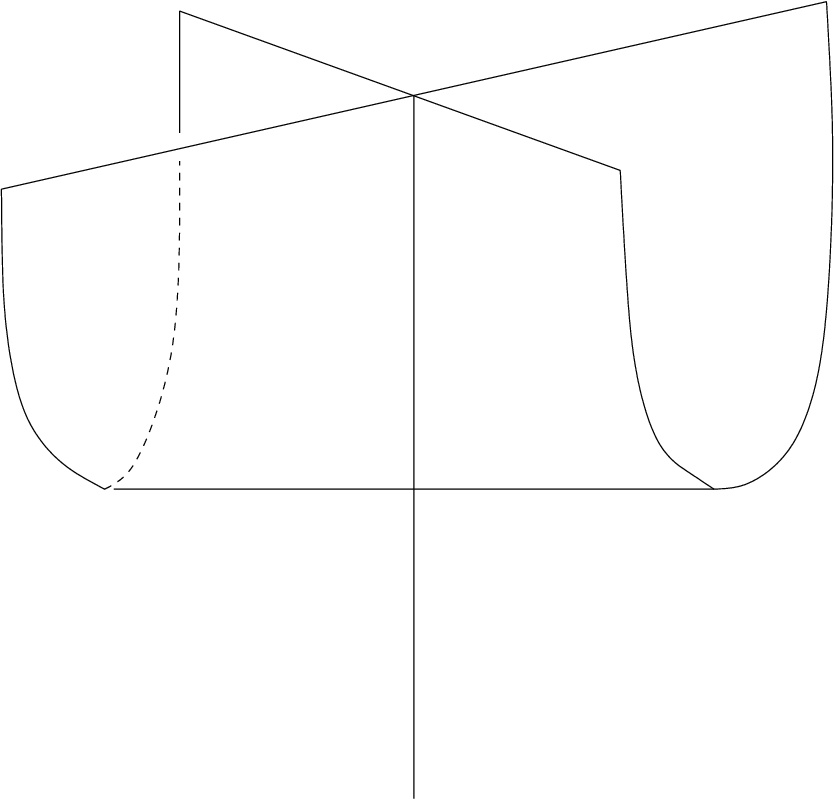, height= 40 mm}
   \caption{Whitney's umbrella}
   \end{figure}

Here again the axis of the variable $z$ 
is contained in the 
algebraic set defined by the equation. The  
half-line where $z<0$ appears separated from 
the points having neighborhoods with  
topological dimension 2, which are precisely those 
verifying the inequality 
$z \geq 0$ in addition of (\ref{Whitneyeq}). 
That is why this
half-line can be  imagined as the  
stick of a (curious) umbrella (see Figure 4). 
We will come back later to
this example, in Subsection \ref{strat}.  
 
Of course, the previous discussion deals with 
phenomena of real
 algebraic geometry.  They do not occur 
in complex geometry. But we
 feel that it is important to develop intuitions 
from visual
 representations of objects, in particular 
from models of real
 algebraic surfaces, even if one 
is mainly interested in complex  ones.

We introduced the cross-cap as a representation 
of the real projective
plane. This representation  
is not faithful, as one identifies like 
this distinct points of the
projective plane. It is a theorem  
of topology that the real projective plane 
cannot be represented
faithfully in space, in the sense  
that it cannot be embedded in $\mathbb{R}^3$: indeed, 
it is non-orientable, and all properly embedded 
surfaces in $\R^3$ are orientable. 
Therefore, if one wants to
represent it in space,  
some singular points are unavoidable, 
in the same way in which
supplementary crossing  
points appear when one represents in the plane 
the abstract graph
which was mentioned at the beginning of this section. 

One faces here a general question, which can be asked 
in any mathematical
category in which a convenient  
notion of \emph{embedding} can be defined: 
\begin{equation}  \label{embed}
    \begin{array}{c}
       \mbox{\emph{Given two objects $X$ and $E$,}} \\ 
       \mbox{\emph{is it possible to represent 
          faithfully $X$ inside $E$? } } 
     \end{array}   
\end{equation}
When the answer is negative, another thing can be 
asked, which leads to
singularity-theoretic  questions:
\begin{equation}  \label{minsing}
    \begin{array}{c}
       \mbox{\emph{Given two objects $X$ and $E$,}} \\ 
        \mbox{\emph{how to represent $X$ 
         inside $E$ with minimal distortion? } } 
     \end{array}   
\end{equation}
For example, Figure 1 shows that the considered 
abstract graph can
be mapped to the plane  
by introducing only one crossing point of the 
simplest type. 
Regarding the projective plane,  
the cross-cap is a more complicated object: 
it is an immersion nearly
everywhere, with the exception  
of the extremities of the segment of 
self-intersection. A
representation could be considered  
simpler if it is entirely an immersion. 
Boy showed in his thesis
\cite{B 03}, done
under the supervision  
of Hilbert, that the projective plane could 
be immersed in
$\mathbb{R}^3$ (see the photographs  
at the end of \cite[VI.48]{HCV 90}).
\medskip

Let us pass now to another category of geometry, 
namely 
\emph{complex algebraic geometry}. A specialization  
of question (\ref{embed}) is: 
\emph{given a smooth complex projective
  curve, is it  
possible to embed it algebraically 
in the complex projective plane?}
We know that this is not always the case, 
as it is  shown by the
following classical theorem: 

\begin{theorem} \label{genusmooth}
   Let $C$ be a smooth algebraic curve inside 
   $\mathbb{CP}^2$. Then
   $g=\frac{(d-1)(d-2)}{2}$, where  
   $g$ is the genus and $d$ is the degree of $C$.
\end{theorem}

This theorem shows that a smooth projective 
curve whose genus 
is not of the form $\frac{(d-1)(d-2)}{2}$ cannot  
be embedded in $\mathbb{CP}^2$.  
Theorem \ref{genusing} below
generalizes this statement to possibly 
singular curves.  

Starting from a smooth curve $C$ 
which cannot be embedded in $\mathbb{CP}^2$,
one can specialize  
question (\ref{minsing}). Here one gets the 
following classical theorem:

\begin{theorem}
    Let $C$ be a smooth projective curve. 
Then there exists an immersion 
    $C\stackrel{\pi}{\rightarrow}\mathbb{CP}^2$ 
whose image has only normal 
    crossings and which is an isomorphism over 
the complement of the
    singular points of $\pi(C)$. 
 \end{theorem}

Here, $\pi(C)$ is said \emph{to have normal 
crossings} if its germ at any singular point  
has  only two irreducible components which 
are smooth and intersect transversely. There is a 
a higher-dimensional version of this notion (see 
Definition \ref{defnormcross}). One may prove the 
theorem by showing that, starting from any embedding 
of $C$ in a projective space, a generic linear 
projection satisfies its conclusions.

What happens if instead of starting from a 
smooth curve, one starts
from a curve which admits  
singular points? Then one has first to ask 
a specialization of
the desire (\ref{fundesire}).  An answer 
to this is the following:

\begin{theorem} \label{resocurve}
    Let $C$ be a projective curve. Then 
there exists a smooth
    projective curve $\tilde{C}$  
    and a morphism $\tilde{C}\rightarrow C$ 
which is an isomorphism
    over the complement  
    of the set of singular points of $C$.
 \end{theorem}

This theorem is historically  the first result 
of \emph{resolution of
  singularities} 
in algebraic geometry. It goes back to the 
construction by Riemann of
the surfaces bearing nowadays his name, 
associated to any algebraic
function of one variable (see the explanations 
which follow
Proposition \ref{vandis}). 
 
 In the sequel, we shall explain different proofs 
of this result, as
 well as of the analogous  
 result for surfaces. But instead of restricting 
to complex projective
 varieties, we shall  
 work with the more general notion of 
\emph{complex analytic spaces} (which we will 
call also shortly \emph{analytic spaces}, or even 
\emph{spaces}). 
 General references about them are e.g. 
the encyclopaedia
 \cite{*** 94}, as well as the books 
of Fischer \cite{F 76} and 
 Kaup \& Kaup \cite{KK 83}. 

 All the spaces we  consider will be assumed 
\emph{reduced}. We will  
 explain everything as intrinsically as possible, 
in order to
 emphasize the various morphisms used in the 
constructions.

\medskip
If $X$ is an analytic space, we  denote by 
$\mbox{Sing}(X)$ its
singular locus. 
 
\medskip
 Section \ref{general} contains the general 
notions necessary 
 to understand the proofs of the theorems 
concerning the existence of 
 resolutions of curves and surfaces 
explained in sections 
 \ref{curve} and \ref{Jung}. In 
Section \ref{open} we state some open
 problems.

\medskip
\section{Generalities about finite morphisms and 
  modifications}
\label{general}

\medskip
\subsection{Finite morphisms} \label{finmor} $\:$
\medskip

When we draw on a piece of paper a real 
surface situated in
3-dimensional space, as we did before for the  
cross-cap and  Whitney's umbrella, 
we trace some curves in the
plane. Let us think for a  
moment about their relation with the surface. 
Suppose that the
drawing is done by cylindrical  
projection to the plane. For the most economic 
drawings, as the one of
Figure 4, one sees that the curves are of three types:

\begin{enumerate}
    \item projections of curves drawn on the surface 
in order to
    cut a part of it; 
    
    \item projections of the curves contained in 
the singular locus of
    the surface; 
    
    \item apparent contours of the surface 
with respect to the chosen
    projection. 
\end{enumerate}

The reader is encouraged to recognize each of 
these types in Figure 4. 

Moreover, there are other important aspects 
of the drawings done before:
each point   
of the plane was the image of 
\emph{ only a finite number of points of the
  surface} and \emph{no point of the surface 
escapes to infinity}. 
This type of  
projection is of great importance in algebraic 
or analytic geometry: 

\begin{definition} \label{definite}
    A morphism $Y \stackrel{\psi}{\rightarrow} X$ 
of reduced complex analytic
    spaces is called \textbf{finite}  if it is 
proper and with finite fibers. 
    
    Let $Y \stackrel{\psi}{\rightarrow} X$ be a finite
    morphism. Suppose moreover that $Y$ is 
equidimensional and that
    $\psi$ is surjective. The \textbf{degree}
    $\mathrm{deg}(\psi)$ of $\psi$ is the 
maximal number of  
    points in its fibers. The \textbf{critical locus}
    $C(\psi)\subset Y$ of $\psi$ is  the set 
of points $p \in Y$ such
    that $\psi$ is not a local analytic 
isomorphism in the
    neighborhood of $p$. The 
\textbf{discriminant locus}
    $\Delta(\psi)\subset X$ 
is the image $\psi(C(\psi))$. 
 \end{definition}

The cardinal of the fibers of $\psi$ 
is equal to $\mathrm{deg}(\psi)$
on the complement of a nowhere dense 
analytic subset of $X$. It is
important to understand that this subset 
is contained in
$\Delta(\psi)$, but that it is not necessarily 
equal to $\Delta(\psi)$. 
Think for example of the normalization morphism 
of an irreducible germ
of curve, a notion explained in the next subsection. 

Already for curves, the notion of discriminant is 
extremely rich,
having a lot of avatars. We recommend Abhyankar's 
fascinating journey
\cite{A 76} among them. We mention also that 
a general program for
studying discriminants 
in singularity theory was described by Teissier 
\cite{T 77} and a general
framework for studying discriminants in 
algebraic geometry was
described by Gelfand, Kapranov \& Zelevinsky 
in \cite{GKZ 94}. 

\medskip

The name ``discriminant locus'' comes from the 
fact that for
projections of hypersurfaces, it is defined by 
the discriminant of a
polynomial :

\begin{proposition} \label{vandis}
    Let $f \in \mathbb{C}[t_1,...,t_{n+1}]$. 
Denote by $Y$ its
    vanishing locus in  
    $\mathbb{C}^{n+1}$, by $X$ the hyperplane 
of $\mathbb{C}^{n+1}$
    defined by $t_{n+1}=0$  and by $\psi$ 
the restriction to $Y$ of the
    canonical projection  
    of  $\mathbb{C}^{n+1}$ onto $X$. 
Then $\psi$ is finite if and only
    if $f$ is  
    unitary with respect to the variable $t_{n+1}$, 
and if this is the
    case, then  
    the discriminant locus of $\psi$ is defined 
by the vanishing of the 
    discriminant of the polynomial $f$ with 
respect to the variable $t_{n+1}$. 
 \end{proposition} 

In the literature one also finds the names
\emph{ramification locus} instead 
of \emph{critical locus} and 
\emph{branch locus} instead 
of \emph{discriminant locus}. 

If $n=1$ in the previous proposition, 
then from the equation $f(t_1,
t_2)=0$ one can express $t_2$ 
as a (\emph{multivalued}) function of
$t_1$. This kind of function 
was called \emph{an algebraic function}
in the XIX-th century. Riemann \cite{R 51} 
associated to such a
function a surface (called 
nowadays \emph{the Riemann surface of the
  function}) over which the function $t_2(t_1)$ becomes
\emph{univalued}. This surface is smooth and 
projects canonically onto
the $t_1$-axis. Riemann explained how one could
construct it by cutting adequately the plane 
along curves connecting
the various points of the discriminant locus, 
which in this case is a
finite set of points on the $t_1$-axis, and by gluing
adequately a finite number of copies of the 
trimmed surface. 
An important point to understand is that 
this Riemann surface does not
project canonically only onto the affine line of 
the independent variable, but also on the affine 
curve of equation
$f(t_1, t_2)=0$, by a map which is a resolution 
of the curve. This is
the reason why we stated in the introduction 
that Theorem
\ref{resocurve} goes back to Riemann.

\medskip
 
 Returning to Definition \ref{definite}, 
the discriminant locus of a
 finite surjective morphism is in fact a 
closed analytic subset of  
 the target space. Moreover, it can be 
naturally endowed with a
 structure of (possibly  
non-reduced) complex space, whose 
formation commutes with base change
 (see Teissier \cite{T 77}). 

 One can construct purely algebraically 
a finite morphism, starting from a 
 convenient sheaf of $\mathcal{O}_X$-algebras. 
In order to
 explain why, recall one  
 of the fundamental ideas of scheme theory: 
an affine algebraic
 variety is \emph{completely  
 determined} as a topological space by its 
algebra of regular
 functions. This gives a  
 procedure to construct spaces by doing algebra: 
each time a new
 algebra (of finite type)  
 is constructed, one gets automatically a new 
affine variety. More
 generally, this can be  
 done over a base which is not an algebraic 
variety, for example over
 a complex  analytic space $X$ 
(see Peternell \& Remmert
 \cite[II.3]{*** 94}). 
In this case, one gets a new complex analytic
 space \emph{over the initial   
 one} $X$ from a quasi-coherent sheaf $\mathcal{A}$ of
 $\mathcal{O}_X$-algebras of finite presentation.  
 The new space is called \emph{the analytic spectrum} 
of the sheaf
 $\mathcal{A}$, and is  
 denoted $\mathrm{Specan}(\mathcal{A})$. 
Denote also by $\pi_{\mathcal{A}}: 
 \mathrm{Specan}(\mathcal{A}) \rightarrow X$ 
the canonical morphism
 associated with  this construction. 
A particular case of it is:

\begin{proposition} \label{finspec}
   If $\mathcal{A}$ is coherent as 
an $\mathcal{O}_X$-module, then the
   morphism $\pi_{\mathcal{A}}$ is finite and
   $(\pi_{\mathcal{A}})_*\mathcal{O}
    _{\mathrm{Specan}(\mathcal{A})}
   \simeq \mathcal{A}$. 
\end{proposition}

Let us state now the property of the 
discriminant loci which relates 
them with the discussion  
about the drawing of surfaces which took place 
at the beginning of
this subsection: 

\begin{proposition}
     Suppose that $Y \stackrel{\psi}{\rightarrow} X$ 
is a  finite surjective
     morphism between equidimensional reduced complex 
analytic spaces
     and that $X$ is smooth. Then its  
discriminant locus is equal,
     set-theoretically, to the union  of the image 
     $\psi(\mathrm{Sing}(Y))$ of the singular 
locus of $Y$ and of the
     closure of the apparent  
     contour (that is, the set of critical values) 
of the restriction
     of $\psi$ to the smooth locus of $Y$. 
\end{proposition}

\medskip
\subsection{The normalization morphism}
\label{normsect} $\:$
\medskip

In the sequel, we shall examine  desire 
(\ref{fundesire}) for
reduced complex analytic spaces.  
We consider that such a space $X$ is ``complicated'' 
if it is
singular. Then we would like  to  
represent it as the image of a non-singular one. 
But we have to
decide first what we want to understand by ``image''. 
The most
encompassing approach would be to consider a  
surjective morphism from any complex analytic space 
onto $X$. But as
we consider that we are happy enough when a point 
is non-singular, it is
natural to ask the morphism to be an isomorphism 
over the set of
smooth points. Such morphisms are particular cases 
of those  which are
isomorphisms over dense open subsets  
(see Peternell \cite[Chapter VII]{*** 94}):

\begin{definition} \label{modif}
    Let $X$ be a reduced complex space. 
    A \textbf{modification} of $X$ is a
    proper surjective morphism   
    $Y \stackrel{\rho}{\rightarrow} X$ such 
that there exists a
    nowhere dense complex subspace  
    $F$ of $X$ with the property:
       $$Y \setminus \rho^{-1}(F) 
  \stackrel{\rho}{\rightarrow} X \setminus F \:
                   \mbox { is an isomorphism.}$$
    The minimal subspace $\mathrm{Fund}(\rho)$ 
with this property is
                   called \textbf{the  
    fundamental locus}  of the modification $\rho$. 
The preimage
                   $\mathrm{Exc}(\rho):=  
    \rho^{-1}(\mathrm{Fund}(\rho))$ 
    of the fundamental locus is called  
\textbf{the exceptional locus}
                   of $\rho$.  
    If $Z$ is a closed irreducible subspace of $X$, 
not contained 
    in the fundamental locus $\mathrm{Fund}(\rho)$, 
then its
                   \textbf{strict transform} 
$Z'_{\rho}$    by the modification 
    $\rho$ is the closure of $\rho^{-1}(Z\setminus
                   \mathrm{Fund}(\rho))$ in $Y$. 
 \end{definition} 

In the literature, strict transforms are also
called \emph{proper  transforms}. 
 
 Informally, to \emph{modify} $X$ means to take 
out a nowhere dense
 analytic subset $F$ and  
 to replace it by another analytic set $E$. 
The important thing to
 remark is that the isomorphism  
 between the ``unmodified'' parts of the 
two spaces must extend to an
 analytic morphism  from the new space to the 
initial one. 

  If $\rho$ is a modification, we can try to 
understand it by looking 
at its fibers, which are {\em compact} analytic spaces
(remember that we asked $\rho$ to be proper!). The
simplest situation arises when all those fibers 
are finite from the
set-theoretical viewpoint, that is, when the 
modification is a finite 
morphism (see Definition \ref{definite}). 
Among such modifications, 
there is a unique one (up to \emph{unique} 
isomorphism) which dominates all
the other ones, the 
\emph{normalization morphism}. 
Before stating precisely this result
(see Theorem \ref{exnorm} below), 
we recall briefly the notion of
\emph{normal analytic space}.

 This concept was first introduced  
in algebraic geometry by Zariski \cite{Z 39}, 
inspired by the
arithmetic notion of integral  
closure and by the notion of 
\emph{normal projective variety} used by
the Italian geometers  
(see Zariski \cite[Footnote 26]{Z 39} and 
Teissier \cite[Section
3.1]{T 95}). It was extended to the complex analytic   
category in the years 1950. Here we prefer to give a
``transcendental''  (that is,  
function-theoretical, non-algebraic) definition, 
which has the
advantage to allow us later on  
to construct very easily holomorphic functions 
on normal varieties. At
the end of the subsection,  
we will briefly come back to Zariski's algebraic 
viewpoint.

The following theorem was proved by Riemann. 
It allows one to show
that a  function of one variable is  
holomorphic on a neighborhood of a point only by 
knowing its behaviour
outside the point.

\begin{theorem} \label{rext}
   (\textbf{Riemann extension theorem})
   Let $U$ be a neighborhood of $0$ in $\mathbb{C}$ 
and $f$ be a
   holomorphic and bounded  
   function  on $U\setminus 0$. Then $f$ extends 
(in a unique way) to
   a function holomorphic on $U$.
 \end{theorem}
 
 The previous theorem, also known as 
\emph{Riemann's removable
 singularity theorem},  was extended to higher 
dimensions
 (see Kaup \& Kaup \cite[chapter 7]{KK 83}):

 \begin{theorem} (\textbf{generalized Riemann 
extension theorem}) \label{rext2}
   Let $U$ be a neighborhood of $0$ in $\mathbb{C}^n, 
   \: n\geq 1$ and
   $f$ be a holomorphic  
   and bounded function  on $U\setminus Z$, 
  where $Z$ is a strict
   closed complex analytic  
   subspace of $U$. Then $f$ extends 
  (in a unique way) to a function
   holomorphic on $U$. 
 \end{theorem}
 
 It is then natural to ask which complex analytic 
 sets admit the same
 property as $\mathbb{C}^n$. In  
 fact, at the beginning of the years 1950, 
 some specialists of complex
 analytic geometry took  
 this property as a definition of a 
 \emph{complex analytic set} (see Remmert
 \cite[pages 30-31]{*** 94}). Later,  
 as this name began to  
 designate any set glued analytically from subsets 
 of $\mathbb{C}^n$
 which are defined locally by a finite  
 number of analytic equations, sets with the Riemann 
 extension property got a special name:

\begin{definition}
    Let $X$ be a reduced complex analytic space. 
    If $U$ is an open
    subspace of $X$,  
    a \textbf{weakly holomorphic function} 
on $U$ is  a holomorphic
    and bounded  
    function defined on $U \setminus Y$, 
where $Y$ is a nowhere dense closed
    subspace of $U$.  
    The space $X$ is called \textbf{normal} 
if every weakly
    holomorphic function  
    on $U$ extends in a unique way to 
a holomorphic function on $U$, and this 
    must occur for any open subset $U$ of $X$. 
 \end{definition}
 
 We have presented the normal spaces as those which 
have in common
 with the smooth ones,  the truth of 
the generalized Riemann
 extension theorem.   In the next theorem 
we state other similarities
 between them: 

 \begin{theorem} \label{propnorm}
      A normal complex analytic space is 
locally irreducible and  smooth in 
      codimension 1 (that is, its singular set 
has codimension $\geq 2$).
  \end{theorem}

 Not any complex analytic space is normal. However, 
any complex analytic
 set can be canonically presented as the image of  
a normal one: 

\begin{theorem} \label{exnorm}
   Let $X$ be a reduced complex space. Then there 
exists a modification 
   $\overline{X} \stackrel{\nu}{\rightarrow}X$ 
such that:
   $\overline{X}$ is normal and $\nu$ is 
a finite morphism. Moreover, 
   if $\nu$ is a fixed modification having 
these properties, then for
   any finite   
   modification $X_1 \stackrel{\nu_1}{\rightarrow}X$, 
there exists a
   unique morphism   
   $\overline{X}\stackrel{\chi}{\rightarrow}X_1$ 
making the following
   diagram commutative: 
           $$\xymatrix{\overline{X}
               \ar[dr]_{\nu} 
               \ar[rr]^-{\chi} & 
               &  X_1
                 \ar[dl]^{\nu_1} \\
               & X & }$$
  \end{theorem}

\begin{definition}
     A morphism $\overline{X} 
\stackrel{\nu}{\rightarrow}X$ as in the previous 
     theorem is called 
\textbf{the normalization morphism} of $X$.
 \end{definition}
 
 Theorem \ref{exnorm} explains why we have used 
the article ``the''
 instead of  
 ``a'' : it implies that a normalization morphism 
is unique up to
 unique isomorphism  above $X$, which is  
 the greatest type of uniqueness in a category. 
In this way, one
 characterizes the normalization morphism by 
a universal property. 

As another consequence of Theorem \ref{exnorm}, 
notice  that the
process of normalization is of \emph{local} nature, 
that is, the
restriction of the 
normalization morphism of $X$ to an open set 
$U\subset X$ is the
normalization morphism of $U$.

 The normalization morphism is a particular 
case of the construction
 of the analytic spectrum 
(see Proposition \ref{finspec}), in which   
 $\mathcal{A}:= \tilde{\mathcal{O}}_{X}$, 
the sheaf of weakly
 holomorphic function  
 on $X$. This sheaf is coherent as an 
$\mathcal{O}_X$-module and can
 be defined algebraically, as was seen already by
 Riemann in the case of complex curves:
 
 \begin{theorem} \label{transcalg}
   Let $X$ be a reduced complex space. 
   The sheaf $\tilde{\mathcal{O}}_{X}$ 
of weakly holomorphic functions
   on $X$ is coherent  
   and equal to the sheaf of integral closures of 
the local rings of
   $\mathcal{O}_X$ in their  total rings of fractions. 
The morphism
   $\pi_{\tilde{O}_X}: \mathrm{Specan}(\tilde{O}_X) 
\rightarrow X$ is
   the normalization morphism of $X$. 
 \end{theorem}

 The \emph{total ring of fractions} $\mathrm{Tot}(A)$ 
of a given ring
 $A$ is by definition the ring 
 of quotients in which all the elements of $A$ which 
are not
 $0$-divisors become invertible. If the initial ring 
is integral, that is,
 without $0$-divisors, then its total ring 
of fractions is a 
 field. If the ring is reduced but not integral, 
that is,
 the associated space is reduced but not irreducible, 
then
 $\mathrm{Tot}(A)$ is canonically the direct product 
of the fields of
 quotients of the rings associated to the irreducible 
components. 
 
 The previous theorem is the key for understanding 
why normalization has 
 at the same time an algebraic and a transcendental 
aspect.

 The concept of normalization is essential when one 
is thinking about 
 resolution of singularities. Indeed, as  shown 
by Theorem
 \ref{propnorm}:

\begin{proposition}
 The  normalization morphism of $X$ 
\emph{separates the local analytically
 irreducible components} of $X$ and 
\emph{resolves the singularities in
 codimension $1$}.  
\end{proposition}

The last statement means that 
$\mathrm{Sing} (\overline{X})$ has codimension 
at least $2$ in $\overline{X}$. 

Let us illustrate the proposition with 
a simple example, that of two
smooth plane curves intersecting transversely, 
met in a topological
context at the beginning of the introduction. 
Here we consider the
union $X$ of the two axes in the complex 
affine plane $\mathbb{C}^2$,
with coordinates $x,y$. Thus, the associated algebra 
is $A:=
\mathbb{C}[x,y]/(xy)$. Consider the function 
$f=x/(x+y)$ 
restricted to $X$. It is weakly holomorphic, 
as it is holomorphic
outside the origin and bounded in a neighborhood 
of it. Theorem
\ref{transcalg} shows that $f$ becomes 
a holomorphic function on the
normalization of $X$. As $f$ is constant outside 
the origin in
restriction to both axes, taking the values $0$ 
and $1$ respectively,
we see that there are two possible limits at 
the origin. Therefore $f$
cannot be extended to a continuous function 
defined all over $X$. The
abstract construction $\mathrm{Specan}(\tilde{O}_X)$ 
separates the
lines, such that $f$ becomes a function holomorphic 
all over the
new curve, which is isomorphic to the disjoint union 
of the two axes. 

To illustrate also Theorem \ref{transcalg}, notice 
that  $f$
is indeed an element of the integral closure of $A$ 
in its total ring
of fractions $\mathrm{Tot}(A)$: $x+y$ is not a 
divisor of $0$ and one
has 
$f^2-f=0$, which is a relation of integral 
dependence of $f$ over $A$.

\medskip

For more details about normal varieties, one may 
consult Greco's book
\cite{G 78}. For details about the more general 
notion of {\em weakly
  normal complex spaces} (in which any 
{\em continuous} weakly
holomorphic function is in fact holomorphic), 
one may consult Adkins,
Andreotti and Leahy's book \cite{AAL 81}. 

 We conclude this subsection with a quotation 
from the introduction
 of Zariski's work \cite{Z 39}: 
 
 \begin{quote}
 
   {\small   Here we introduce the concept 
of a  \emph{normal}
 variety, both in the affine and in the  
 projective space, and we are led to a geometric 
interpretation of the
 operation of integral  
 closure. The importance of normal varieties is 
due to: ... \emph{the
 singular manifold of a normal $V_r$ is of dimension 
 $\leq r-2$ } (in particular a normal  
 curve $(V_1)$ is free from singularities)... 
There is a definite
 class of normal varieties associated with and  
 birationally equivalent to a given variety $V_r$. 
This class is
 obtained by a process of  
 integral closure carried out in a suitable fashion 
for varieties in
 projective spaces....  
 
 The special birational transformations effected 
by the operation of
 integral closure, and the  
 properties of normal surfaces, play an essential 
 r{\^o}le in our
 arithmetic proof for the reduction of  
 singularities of an algebraic surface.  }
 \end{quote}

\medskip
\subsection{Blowing-up points and subschemes}  
\label{blow}$\:$
\medskip

Let us begin by a theorem of elementary 
geometry (see Figure 5):

\begin{proposition}
      Let $ABC$ be a triangle in the euclidean plane. 
For each point
      $P$ in the plane, 
      consider the symmetric lines of $PA, PB, PC$ 
with respect to
      the bisectors  
      of the angles $\angle BAC, \angle CBA$ 
and $\angle
      ACB$ respectively. Then these three new lines   
      intersect at another point $s(P)$ and the 
transformation $P
      \rightarrow s(P)$ is an involution.  
 \end{proposition}

 The proposition can be easily proved using the 
classical theorem of
 Ceva. It is also true that  
 $P$ and $s(P)$ are the two foci of a conic tangent 
to the edges
 of the triangle $ABC$ (as an  
 illustration of this fact, think at the inscribed 
circle, which is a
 conic tangent to the  
 three edges, and whose center $I$ verifies 
$I = s(I)$). 
 
 But what interests us here more is the fact 
that \emph{the mapping $s$
 is not defined everywhere}.  
 Indeed, it is not defined at the vertices of 
the triangle. By doing
 some drawings, one sees  
 experimentally why:  
 if one tends to a vertex by remaining on a line 
passing through it, then the limit of the  
 transforms is well-defined, but 
\textit{it depends on the chosen
 line}. Moreover, by varying the line,  
 one gets as limits all the points situated on 
the line containing the opposite
 edge. Therefore, in a way:   
 $$s \: \: \mbox{\emph{transforms each vertex 
into the opposite edge.}}$$
 As the
 dimension increases like this from $0$ to $1$, 
one assists to a kind of 
 ``blowing-up'' of each vertex into a line. 
One has at the same time a phenomenon 
of ``blowing-down'' of each 
edge of the triangle into the opposite vertex.
Indeed, all the points of the line containing 
an edge, with the only exception of the 
vertices, are sent by $s$ into the opposite 
vertex of the triangle.  

This kind  of examples 
led Zariski to introduce a general notion of 
``blowing-up'' and ``blowing-down'' in algebraic 
geometry. In order to
 explain it, let us 
 first express algebraically the transformation $s$.

 There are other points $P$ for which 
$s(P)$ is not defined,
 those  for which the  
 three new lines are parallel. 
But in this case $s(P)$ can be
 interpreted as a point at  
 infinity, which shows that it is better 
to think about $s$ as a
 transformation of the \emph{projective}  
 plane into itself. There is then a choice 
of coordinates
 which makes the  
 transformation particularly simple from 
the algebraical viewpoint:
 choose the  unique system of 
 projective coordinates $(X:Y:Z)$ such 
that the equations $X=0, Y=0,
 Z=0$ define the  edges of the  
 triangle, and such that the center $I$ of 
the incircle is $(1:1:1)$. Then the
 involution $s$ can be written as:
 \begin{equation} \label{quadreq}
   (X:Y:Z) \cdots\rightarrow (\frac{1}{X}: 
   \frac{1}{Y}: \frac{1}{Z}).
 \end{equation}

{\tt    \setlength{\unitlength}{0.92pt}}
  \begin{figure}
  \epsfig{file=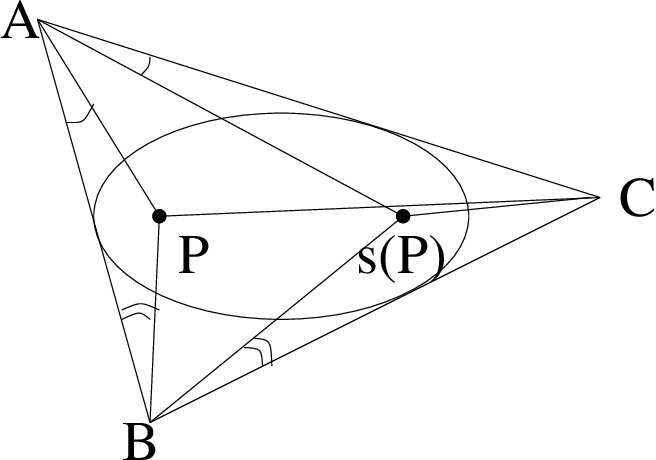, height= 40 mm}
  \caption{A birational involution of the plane}
 \end{figure}

Since the same map can be expressed as 
$(X:Y:Z) \cdots\rightarrow
(YZ: ZX: XY)$, that is, as a map 
with quadratic polynomials as
coordinates, one speaks about a 
\emph{quadratic transformation of}
$\mathbb{P}^2$. For a deeper understanding of  
this vocabulary, we
refer the reader to the quotation from 
Zariski \cite{Z 42} at the end of this subsection. 

We see that $s$ can be expressed in projective 
coordinates using
rational functions of the coordinates. 
That is why one says that $s$ is a
\emph{rational map}. Because its inverse is 
also rational (as the map $s$ is an
involution), one says that the map 
$\mathbb{P}^2\stackrel{s}{\cdots
  \rightarrow}\mathbb{P}^2$ is 
\emph{birational}. Generally speaking:

\begin{definition}
    Let $X$ and $Y$ be two reduced and irreducible 
algebraic varieties. 
    A \textbf{rational map} 
$Y\stackrel{s}{\cdots\rightarrow}X$ is an 
    algebraic morphism $U \rightarrow X$, 
where $U$ is a dense Zariski
    open set of $Y$. \textbf{The indeterminacy locus} 
of a rational
map is the complement of the maximal possible such $U$.

    A \textbf{birational map}
    $Y\stackrel{s}{\cdots\rightarrow}X$ is 
a rational map which
    realizes an isomorphism between dense 
open subsets of $Y$ and $X$. 
    A \textbf{birational morphism} 
is a birational map which is  defined everywhere.
\end{definition}

\emph{Birational geometry} is the study of 
algebraic varieties up to
birational isomorphism. It seems to have begun 
as a conscious 
domain of research with Riemann's definition 
\cite[chapter XII]{R 57} 
of the birational
equivalence of plane algebraic curves, which 
we quote here:

 \begin{quote}
 
   {\small  We shall consider now, as pertaining 
to a \emph{same
   class, all the irreducible algebraic equations 
between two variable
   magnitudes, which can be transformed the ones 
into the others by
   rational substitutions.}}

 \end{quote}

The notion of \emph{modification} (see Definition 
\ref{modif}) was
introduced in complex analytic  
geometry in order to extend to it the notion 
of birational
morphism, and to create an analog of 
the birational geometry, the
so-called \emph{bimeromorphic geometry}. 
To understand this, notice
that a proper birational morphism $Y \rightarrow X$  
between complex algebraic varieties is a 
modification of the underlying complex analytic 
space of $X$.

By definition, the difference between 
the concepts of \textit{rational
  map} and \textit{rational 
morphism}  is that for the first one we allow 
the presence of points 
of indeterminacy, while this is forbidden for 
the second notion. There is a
general way to express a rational map in terms 
of rational
morphisms. One simply 
considers \emph{the closure of the graph of 
the rational map}. As this
closure lives in the product space, it can 
be naturally projected 
onto the factor spaces, which are the source 
and the target of the
initial map. But these two projections are 
now \emph{morphisms}: 

\begin{equation} \label{elimindet}
     \xymatrix{
         & \overline{\mbox{Graph}(s)}
        \ar[dl]_{p_Y} \ar[dr]^{p_X} & 
        \hookrightarrow  Y \times X \\
       Y & \stackrel{s}{\cdots \rightarrow} & X }
\end{equation}

The first one $ \overline{\mbox{Graph}(s)} 
\stackrel{p_Y}{\longrightarrow}
Y$ is a birational \emph{morphism} and the 
second one is also a morphism, but not necessarily 
birational. The map
$s$ can be expressed as the  
composition $s= p_X \circ p_Y^{-1}$. If $X$ is 
complete (that is, its
underlying analytic space is compact), 
then $p_Y$ is proper, and
therefore $p_Y$ is a \emph{modification} 
of the underlying analytic
space of $Y$ (see Definition \ref{modif}).

As a very important example, let us consider 
the canonical projection
map from a vector space $V$ of dimension $n\geq 2$ 
to its
projectivization $\mathbb{P}(V)$: 

\begin{equation} \label{blowpt}
\xymatrix{& B_0(V)
        \ar[dl]_{p_V} \ar[dr]^{p_{\mathbb{P}(V)}} & 
           \hookrightarrow  V \times \mathbb{P}(V)  \\
       V & \stackrel{s}{\cdots \rightarrow} &  
                   \mathbb{P}(V)}   
\end{equation}

One can study this diagram using a fixed coordinate 
system. Start from a basis  
of $V$, which determines an isomorphism between $V$ 
and $\mathbb{C}^n$, and the associated 
canonical covering of $\mathbb{P}(V)$ with $n$ 
affine charts
isomorphic to $\mathbb{C}^{n-1}$. This gives 
a covering of $V \times
\mathbb{P}(V)$ with $n$ charts isomorphic to
$\mathbb{C}^{2n-1}$. Being the
roles of the different coordinates completely 
symmetric, one sees
that it is enough to study the modification 
$p_V$ inside one of these
charts. One proves in this way:

\begin{proposition} \label{blowpoint}
   1)The algebraic variety $B_0(V)$ is smooth. 

   2)The indeterminacy locus
   of the modification $p_V$ is the point $0$ 
and its exceptional
   locus is sent isomorphically to $\mathbb{P}(V)$ 
by the morphism 
   $p_{\mathbb{P}(V)}$. Moreover, this second morphism 
is canonically
   isomorphic to the projection map of the total 
space of the
   tautological line bundle 
$\mathcal{O}_{\mathbb{P}(V)}(-1)$.

   3) If $y_i := \frac{x_i}{x_n}, \: \forall \: i
   \in \{1,..., n-1\}$ are the coordinates of the 
canonical
   chart $U_n := \mathbb{P}(V)\setminus \{x_{n}=0\}$  of
   $\mathbb{P}(V)$, then the canonical projection 
of the affine space
   $V\times U_n$ with coordinates 
$x_1,...,x_n,y_1,...,y_{n-1}$ onto
   the space with coordinates $y_1,...,y_{n-1},x_n$ 
is an isomorphism
   when restricted to $B_0(V)$.

   4) In terms of the coordinates 
$x_1,...,x_n,y_1,...,y_{n-1}$ of 
      $B_0(V)\cap (V \times U_n)$ and $x_1,...,x_n$ 
of $V$, the
      modification  $p_V$ is expressed as:
     \begin{equation} \label{blovect}   
          x_1=y_1 \cdot x_n\: , \: ... \: ,\:  
       x_{n-1}= y_n \cdot x_n\: , \:
          x_n=x_n.
     \end{equation}
\end{proposition} 

This proposition shows that one has modified $V$ 
by replacing the
origin with the projective space of all the 
directions of lines
passing through the origin. Therefore, 
the origin has been ``blown-up''
into a higher dimensional space:

\begin{definition}
  The birational morphism 
$B_0(V)\stackrel{p_V}{\rightarrow} V$ of 
diagram (\ref{blowpt}) is
  called \textbf{the blowing-up of the origin in} $V$. 
\end{definition}

In Figure 6 we have represented the blowing-up 
of the origin in a real
plane, by drawing its restriction over a disc 
centered at the
origin. It is an excellent exercise to 
understand why one gets like
this a M{\"o}bius band. We have represented also 
the strict transforms
$L_i'$ of four segments $L_i$ passing through 
the origin. Please
contemplate how they become disjoint on 
the blown-up disc!

{\tt    \setlength{\unitlength}{0.92pt}}
  \begin{figure}
  \epsfig{file=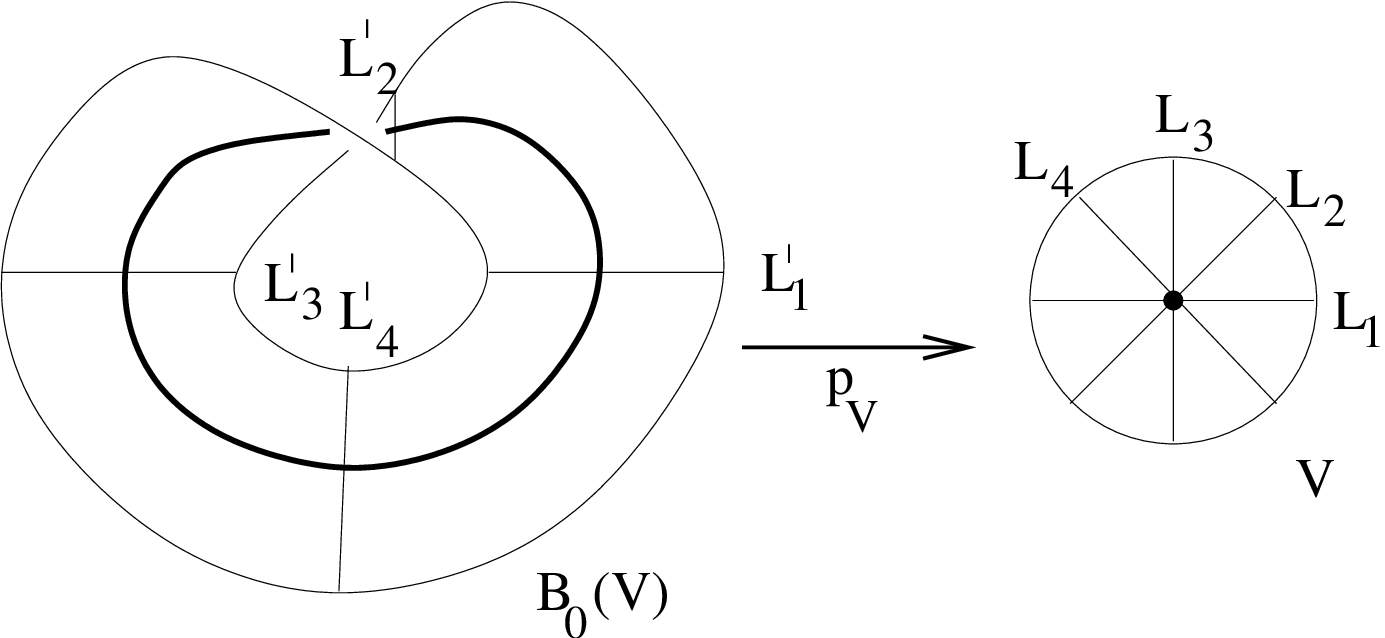, height= 40 mm}
  \caption{Blowing-up the origin in a real plane}
   \end{figure}

The construction of the blowing-up of a point may 
be extended from an
ambient vector space to  
an arbitrary complex manifold. One may blow-up 
a point of it by
choosing a system of local coordinates and 
by identifying like this
the point with the  origin of the vector space 
defined by that
coordinate system. Different coordinate systems 
lead to blown-up
spaces which are canonically isomorphic over 
the initial manifold,
which shows that the blow-up exists and 
is unique up to unique
isomorphism.

Roughly speaking, one blows-up a point of 
a smooth surface by replacing it
with a rational 
curve whose points correspond to the projectified 
algebraic tangent
plane of the surface at that point. In the same way, 
one blows-up a
point in a complex manifold by replacing it with 
the projectified
tangent space at that point. More generally, 
one can blow-up a
submanifold by replacing it with 
its projectified normal bundle. But
one can still generalize this construction, 
and blow-up a
non-necessarily smooth and even non-necessarily 
reduced subspace. The
following theorem, characterizing blowing-ups by 
a universal property, was proved by Hironaka 
\cite{H 64}:

\begin{theorem}
    Let $X$ be a (not necessarily reduced) 
complex analytic space. 
    Let $Y$ be a subspace of $X$, defined by the  
    ideal sheaf $\mathcal{I}$. Then there exists 
a modification 
    $B_Y(X) \stackrel{\beta_{X,Y}}{\rightarrow} X$ 
such that:
    
    $\bullet$ the preimage ideal sheaf 
$\beta_{X,Y}^{-1}\mathcal{I}$
    is locally invertible;
    
    $\bullet$ for any morphism 
$B\stackrel{\beta}{\rightarrow}{X}$
              such that  
              $\beta^{-1}\mathcal{I}$ 
is locally invertible, there
              exists a unique morphism  
              $\gamma$ such that the following 
diagram is commutative:
       $$\xymatrix{B
         \ar[dr]_{\beta} 
             \ar[rr]^-{\gamma} & 
       &  B_Y(X) 
                 \ar[dl]^{\beta_{X,Y}} \\
    & X & }$$
\end{theorem}

\begin{definition}
    A modification $B_Y(X) \stackrel{\beta_{X,Y}}
{\rightarrow} X$ as in
    the previous theorem   
    is called \textbf{the blowing-up of $Y$} 
(or \textbf{with center
    $Y$}, or \textbf{of $\mathcal{I}$) in $X$}. 
\end{definition}

In algebraic geometry, blowing-ups are also 
known as \emph{monoidal
  transforms} (see the quotation from Zariski 
at the end of this 
  subsection) and in analytic geometry as
  \emph{$\sigma$-processes}. 

Not all the modifications can be obtained 
by blowing-up a
subspace. Those which can are precisely the 
\emph{projective
  modifications}. Moreover, a blowing-up 
does not determine the
ideal sheaf $\mathcal{I}$ used to define it. 
In fact, a blowing-up in
the sense of Hironaka \cite{H 64} is the couple $(Y,
\beta_{X,Y})$. Notice that the ideal sheaf $\mathcal{I}$
giving birth to it has been forgotten.

Let us consider again the example of the birational 
involution 
(\ref{quadreq}). We have seen that its indeterminacy 
locus in the
projective plane is the set of 
vertices of the initial triangle. Moreover, we saw 
that the
indeterminacy was caused by the fact that when one 
tends to a vertex along
different lines passing through the vertex, 
one gets different limits
of their images by the involution. This suggests  
that, by 
replacing each  vertex with a curve parametrizing the 
lines passing
through it, that is by blowing-up the three vertices, 
one modifies the
projective plane in such a way that now 
the rational map is defined everywhere. 
This is indeed the case, as
shown by the following quotations from 
the article \cite{Z 42} in which Zariski
introduced the operation of blowing-up under 
the name of ``monoidal
transformation'':

\begin{quote}
    {\small With some non-essential modifications, 
and without their
    projective trimmings, the space Cremona 
transformations, known as
    monoidal transformations, are monoidal 
transformations in our
    sense.... A \emph{quadratic transformation} 
is a special case
    of a monoidal transformation, the center 
is in that case a
    point.... A quadratic Cremona transformation 
is not at all a
    quadratic transformation in our sense. 
Our quadratic
    transformation has only one ordinary 
fundamental point and its
    inverse has no fundamental points at all, 
while a plane quadratic
    transformation and its inverse both have 
three fundamental points,
    which in special cases may be infinitely 
near points.... Of
    course, an ordinary quadratic 
transformation between two planes
    $\pi$ and $\pi'$ can be expressed 
as a product of quadratic
    transformations in our sense, 
or more precisely as the product of
    $3$ successive quadratic transformations 
and of $3$ inverses of
    quadratic transformations.}
\end{quote}

\medskip
\subsection{The definitions of resolution and 
embedded resolution}$\:$
\medskip

Coming back to the desire (\ref{fundesire}), 
one could think that it
is fulfilled for complex analytic spaces 
if one finds a
modification whose initial space is smooth. 
Such modifications have a
special name:

\begin{definition} \label{resol}
     Let $X$ be a reduced space. 
A \textbf{resolution} (or
     \textbf{resolution of the  
     singularities}) of $X$ is a morphism
     $\tilde{X}\stackrel{\pi}{\rightarrow} X$ 
such that: 
     
      $\bullet$ $\pi$ is a modification of $X$;
      
      $\bullet$ $\tilde{X}$ is smooth;
      
      $\bullet$ the restriction 
$\tilde{X}\setminus \pi^{-1}(\mathrm{Sing}(X))
                     \stackrel{\pi}{\rightarrow} 
X\setminus
                     \mathrm{Sing}(X)$ 
is an isomorphism. 
 \end{definition}

In some cases, people skip the last condition 
in the definition of a
\emph{resolution}.  Other terms  
which were used in the literature are 
``reduction of the
singularities'' (see Walker \cite{W 35} 
or the quotation at the end of
Subsection \ref{normsect}) and ``desingularization''.

In the previous definition, if $X$ is embedded 
in a smooth ambient
space, one needs sometimes to get a resolution 
as a restriction of a
modification of the ambient space. 
It is also important  to have a
modification in which the subspaces of interest  
are as simple as possible from a local viewpoint. 
The subspaces we are 
speaking about are the exceptional locus of 
the modification and the
strict transform of the space $X$ 
(see Definition \ref{modif}). As the
ambient space is supposed 
to be smooth, it can be shown that the exceptional 
locus is
necessarily of codimension 1. 
The usual condition of local simplicity
for a hypersurface is that of being a 
\emph{divisor with normal
  crossings}, and for the union of a 
hypersurface and a subvariety of
possibly lower dimension, of \emph{crossing normally}:

\begin{definition} \label{defnormcross}
  Let $M$ be a complex manifold and $D$  
a divisor of $M$. Denote 
  by $|D|$ the underlying reduced hypersurface of $D$. 
  
  One says that 
  $D$ is a \textbf{divisor with normal crossings} 
if for each point
  $p\in |D|$,  
  there is a system of local coordinates 
centered at $p$ such that in some 
  neighborhood of $p$ the hypersurface $|D|$ 
is the union of some hyperplanes 
  of coordinates. 
  
  If $V$ is a reduced subvariety of $M$ of 
a possibly lower 
  dimension, one says that $D$ 
\textbf{crosses $V$ normally}   if 
  for each point $p\in |D|\cap V$, 
  there is a system of local coordinates centered 
at $p$ such that in some 
  neighborhood of $p$ the hypersurface $|D|$ is 
the union of some hyperplanes 
  of coordinates and $V$ is the intersection of 
some of the remaining 
  hyperplanes of coordinates. 
\end{definition}

In order to understand this definition, remark that 
if one takes as
the hypersurface  
$|D|$ the union of two of the coordinate hyperplanes 
in $\mathbb{C}^3$, 
then no line intersecting them only at the origin 
crosses them
normally. But after  
blowing-up the intersection of the two planes, 
the total transform of
their union and the strict transform of 
the line cross normally. 

In the previous definition, we do not suppose 
that each irreducible 
component of $D$ is smooth. If this is moreover 
the case, one says 
usually that $D$ is a \emph{strict normal 
crossings divisor}.

\begin{definition} \label{embresol}
     Let $X$ be a closed reduced subspace of 
a complex manifold $M$. An
     \textbf{embedded  
     resolution} of $X$ in $M$ is a morphism
     $\tilde{M}\stackrel{\pi}{\rightarrow} M$ 
such that: 
     
      $\bullet$ $\pi$ is a modification of $M$;
      
      $\bullet$ $\tilde{M}$ is smooth; 
      
      $\bullet$ the restriction 
$\tilde{M}\setminus \pi^{-1}(\mathrm{Sing}(X))
                     \stackrel{\pi}{\rightarrow} 
M\setminus
                     \mathrm{Sing}(X)$ is an 
isomorphism; 
                     
      $\bullet$ the strict transform $X_{\pi}'$ 
is smooth;
      
      $\bullet$ the exceptional locus of $\pi$ 
has normal crossings
          and also crosses  
          normally the strict transform   
  $X_{\pi}'$ of $X$ by the modification $\pi$.
 \end{definition}

\medskip
\subsection{The special case of surfaces} 
   \label{spesurf} $\:$
\medskip

From now on, we will restrict our considerations 
to surfaces. Inside them, it will be important 
to consider also various curves and their 
intersection numbers. That is why we recall 
first basic facts about intersection theory 
of curves on smooth surfaces. For more details, 
one may consult Hartshorne's book  
\cite[Chapter V.1]{H 77}.

Let $C_1$ and $C_2$ be two (not necessarily reduced) 
properly 
embedded curves (that is, effective divisors) 
contained 
on a (possibly non-compact) 
smooth surface $M$. Let $P$ be a common 
point of $C_1$ and $C_2$. Suppose that their  
analytic germs at $P$ have no common components. 
Then their \emph{intersection number at $P$} 
may be defined by the formula:
  $$ (C_1 \cdot C_2)_P := \dim_{\C}
        \mathcal{O}_{M,P}/ (f_1, f_2) >0, $$
where $f_i\in \mathcal{O}_{M,P}$ denotes a 
holomorphic function defining the curve $C_i$ in 
a neighborhood of $P$ on $M$. This definition 
is independent of the choice of the functions 
$f_1, f_2$. 

If $C_1$ is compact (but not necessarily $C_2$), 
then one may define the \emph{global intersection 
number}:
  $$ (C_1 \cdot C_2) := \deg \mathcal{O}_M(C_2) 
       |_{C_1},$$
where $\mathcal{O}_M(C_2)$ denotes the line bundle 
generated by the divisor $C_2$ on $M$ and 
$\mathcal{O}_M(C_2)|_{C_1}$ denotes its 
restriction to $C_1$. This number depends 
only on the germ of $C_2$ along $C_1$. 
If \emph{both} $C_1$ and $C_2$ are compact, 
then one has :
  $(C_1 \cdot C_2) = (C_2 \cdot C_1)$. 

By bilinearity, one may extend the definition of 
intersection number to the case where 
$C_1$ and $C_2$ are possibly non-effective 
compact divisors on the surface $M$.

If 
$C_1$ and $C_2$ have no common components, then 
they have only a finite set of common points and:
 $$(C_1 \cdot C_2) = \sum_{P \in C_1 \cap C_2} 
      (C_1 \cdot C_2)_P. $$
This shows that, with the previous hypothesis, 
$(C_1 \cdot C_2) >0$ whenever 
$C_1 \cap C_2$ is non-empty. 

This positivity result is no longer necessarily 
true if $C_1$ and 
$C_2$ have common components, for 
example if $C_1 =C_2$. The simplest example 
is provided by the self-intersection number 
of the exceptional curve created by blowing-up a 
point on a smooth surface, which is equal to $-1$ 
(see below). 

\medskip

If $M$ is a 
smooth complex surface and $D$ is a compact 
reduced curve in $M$, 
one associates to it an abstract unoriented 
\emph{weighted dual graph} $\Gamma(D)$. 
Its vertices $v_i$ correspond  
bijectively to the irreducible components $D_i$ 
of $D$ and 
its edges correspond bijectively to 
the unordered pairs of distinct vertices $\{v_i, v_j\}$ 
whose corresponding 
irreducible components $D_i, D_j$ 
intersect. Each vertex $v_i$ is weighted 
by the self-intersection number $e_i := D_i^2$
of the corresponding component and each 
edge by the intersection number of the 
components associated to its vertices.
Denote by $e_{i,j}=e_{j,i}:= D_i \cdot D_j$ 
the weight of the edge joining $v_i$ and $v_j$. 
Therefore, if there is no edge between 
$v_i$ and $v_j$, one has $e_{i,j}=0$. 

One may associate to the curve $D$ 
\emph{the intersection form} on the free 
abelian group of the divisors supported by $D$, 
given by the intersection number. This 
intersection form depends only on the 
associated dual graph $\Gamma(D)$. Indeed, 
if $\sum_{i} x_i D_i$ is a divisor supported by $D$, 
then its self-intersection number is :
  $$ \sum_{i} e_i x_i^2 + 2 
         \sum_{\{i,j| \: i \neq j\} }e_{i,j}x_i x_j.$$

One can particularize the previous constructions 
to the case where $D$ 
is the exceptional divisor of a resolution 
of a normal surface singularity. 
The divisors appearing like this are very special, 
as is shown by the 
following theorem. Point 1) was proved by 
Du Val \cite{V 44} and 
Mumford \cite{M 61}, 
and point 2) was proved by 
Grauert \cite{G 62}.

\begin{theorem} \label{charex}
    1) Let $D$ be the exceptional locus of 
a resolution of a complex 
         analytic normal surface singularity. 
Then $D$ is a connected
         curve whose intersection form is negative 
definite.
         
     2) Let $D$ be a reduced divisor with compact 
support in a smooth 
          complex analytic surface. If 
the intersection form of $D$ is 
          negative definite, then there exists 
a neighborhood of $D$ which is 
          the resolution of a normal surface 
having only one singular 
          point, such that $D$ is the 
exceptional divisor of this resolution. 
 \end{theorem} 

When the hypothesis of point 2) are satisfied, 
one says that $D$ can
be \emph{contracted}.

For more details about dual graphs and intersection 
matrices, one can
consult Laufer \cite{L 71}, N{\'e}methi \cite{N 99} 
or Popescu-Pampu
\cite{PP 06}. 

Let us consider again the blowing-up of a point 
on a smooth surface,
which is illustrated in Figure 6 for the case 
of real surfaces. As a
particular case of point 2) of 
Proposition \ref{blowpoint}, one shows that 
its exceptional locus $E$ 
is a smooth rational curve of self-intersection 
number $E^2=-1$. Such
a curve passing only through smooth points of 
a surface is called
classically an \emph{exceptional curve of 
the first kind}. More
generally, an \emph{exceptional curve} 
is a reduced divisor which can be
contracted.  It can be shown that 
an exceptional curve of the first
kind must contract to a smooth point of the new normal
surface. 

Moreover, it is a classical theorem of Castelnuovo 
that if
one starts from a (not necessarily smooth) 
projective surface, then
the surface obtained after having contracted 
an exceptional curve of
the first kind is again projective. This is 
to be contrasted with the
general case, when the contraction of an 
exceptional curve cannot be 
always done in the  projective category, 
or even in the category of 
schemes (see an example of Nagata in 
B\u adescu \cite[chapter 3]{B 01}). 
\medskip

Let $(S,s)$ be a normal surface singularity. 
In Section \ref{Jung} we
explain a proof of Theorem \ref{resurf}, 
which says in particular that
a resolution of $(S,s)$ always exists. 
It is then natural to try to
compare all possible resolutions. We have 
the following theorem
concerning them:

\begin{theorem} \label{minresurf}
  Let $(S, s)$ be a germ of normal surface. 
  There exists a minimal resolution $S_{min}
  \stackrel{\pi_{min}}{\longrightarrow}S$  
of $(S,s)$, in the sense
  that any other resolution 
$S'\stackrel{\pi'}{\rightarrow}S$ can be
  factored through a composition $\gamma$ 
of blowing-ups of points: 
   $$\xymatrix{S'
         \ar[dr]_{\pi'} 
             \ar[rr]^-{\gamma} & 
       &  S_{min}
                 \ar[dl]^{\pi_{min}} \\
    & S & }$$
   The minimal resolution can be characterized 
by the property that no
         irreducible component of $E_{min}$ 
is exceptional of the
         first kind.
\end{theorem}

The previous theorem is specific to surfaces: it is 
no longer true in higher dimensions.

\medskip
\section{Resolutions of curve singularities} 
\label{curve}

\medskip
\subsection{Abstract resolution}$\:$
\medskip

Using theorems \ref{propnorm} and \ref{exnorm}, 
we get immediately:

\begin{theorem}
   If $C$ is a reduced analytic curve, 
then its normalization morphism is a 
   resolution of $C$.
 \end{theorem}
 
 Analytically, a normalization  of a germ 
of curve $(C,c)$  is given by 
 a set of parametrizations 
 $(\mathbb{C},0) \stackrel{\nu_i}{\rightarrow} C_i$ 
of the irreducible 
 components of $C$, with the condition that each 
parametrization 
 realizes a homeomorphism onto its image. If 
an irreducible 
 germ of curve is embedded  in some space 
$(\mathbb{C}^n,0)$, such a 
 parametrization is given by $n$ convergent 
power series 
 $x_1(t),..., x_n(t)$ in a variable $t$, with 
the restriction that one cannot 
 write them as convergent power series of 
a new variable $w$, with 
 $w$ a convergent power series of $t$ 
of order $\geq 2$.  
 
 Let us explain how one can deduce the existence of 
a normalization 
 of an analytically irreducible germ of curve 
with topological arguments,
 in the spirit  
 of Riemann. Consider an embedding of the germ in 
a smooth space, 
 and choose local coordinates $(x_1,...,x_n)$ in 
this space such that 
 the canonical projection onto the axis of 
the first coordinate $x_1$ 
 is finite (which means that the curve is not 
contained in the
 hyperplane of the  
 other coordinates). Denote by
 $(C,c)\stackrel{\alpha}{\rightarrow}(\mathbb{C},0)$  
 the restriction of this projection to the curve. 
Look at the induced morphism 
 of fundamental groups  
 $\pi_1(U\setminus c)
 \stackrel{\alpha_*}{\longrightarrow}
\pi_1(V\setminus 0)$,  
 where $U$ and $V$ are neighborhoods of $c$ in 
a representative of $C$ 
and of $0$ in $\mathbb{C}$,  
 which are homeomorphic to discs. Since the covering 
is finite and has a
 connected   total space, the image group  
 $\alpha_*(\pi_1(U \setminus c) ) \subset 
\pi_1(V\setminus 0)$ is infinite 
 cyclic and has a  finite index $m\geq 1$. 
Notice that $m$ can also be
 interpreted as the degree of the finite 
morphism $\alpha$. 
 
 Consider then another copy of $\mathbb{C}$, 
with parameter $t$, and
 the morphism  
 $\mathbb{C}_t\stackrel{\tau}{\rightarrow}
\mathbb{C}_{x_1}$ defined by the
 equation $x_1=t^m$.  By construction, 
 $\alpha_*(\pi_1(U \setminus c) ) = 
\tau_*(\pi_1(\mathbb{C} \setminus 0) )$, 
 which shows that the map $\tau$ can be lifted 
to a homeomorphism $\nu$ from 
 a pointed neighborhood of the origin in 
$\mathbb{C}_t$ to $U\setminus
 c$. Compose this morphism $\nu$ with the  
 ambient coordinate functions $(x_1,...,x_n)$ 
at the target. By construction, 
 all the functions $x_i\circ\nu$ are holomorphic 
and bounded on a 
 pointed neighborhood of $0$. By Riemann's extension 
theorem 
 \ref{rext2}, all of them can be extended 
to functions holomorphic also at the 
 origin. This shows that $\nu$ extends to 
a map holomorphic all over 
 the chosen neighborhood of $0$ in $\mathbb{C}_t$:
 
 \begin{equation} \label{ramnorm}
         \xymatrix{
           (\mathbb{C}_t,0) \ar[r]^{\nu} 
\ar[dr]_{\tau}& (C,c)
           \ar[d]^{\alpha}\\ 
             &  (\mathbb{C}_{x_1},0) }
  \end{equation}
  
  The map $\nu$ constructed in this way is 
the normalization morphism of
  $(C,c)$. 
\medskip

The problem with this resolution process by 
the normalization morphism is that, 
given a fixed embedding of the curve, 
it does not extend 
naturally to a modification of the ambient space. 
But in many applications, 
and in particular for
Jung's method of resolution of surfaces presented 
in Section
\ref{Jung}, it is important 
to resolve the curve by a morphism which 
is the restriction of an
ambient one. 

Let us consider this second problem in the case 
of plane curves. 
At first, Max Noether
\emph{simplified} the singularities of plane curves 
by doing sequences
of quadratic transforms of the type (\ref{quadreq}), 
with respect to
conveniently 
chosen triangles (see \cite{N 73} and 
the obituary \cite{CES 25} by
Castelnuovo, 
Enriques \& Severi).  If we used the 
term ``simplified'' and not
``resolved'' in the previous sentence, it is because 
he did not really
resolve them 
with the modern definition \ref{resol}. 
He proved instead the theorem:

\begin{theorem} \label{quadr}
  Let $C$ be a plane curve. Then one can transform 
the curve $C$ into
  another curve $C'$ which has only ordinary 
multiple points, by a
  sequence of involutions isomorphic to 
the involution (\ref{quadreq}).
\end{theorem}

An \emph{ordinary multiple point} designates 
a point of the curve at
which its analytically irreducible components 
are smooth and pairwise
transverse.  The strategy to prove the previous 
theorem was to iterate 
the following steps, given the curve 
$C\subset \mathbb{P}^2$ 
we want to simplify:

1) Choose a singular point $c$ of $C$ which is not 
an ordinary 
      multiple point. 
   
2) Choose a triangle having a vertex at $c$, and whose 
edges 
      are transverse to $C$ (that is, they cross 
it normally) outside 
      the set of its vertices.
      
3) Choose a quadratic transformation  of the plane 
whose 
      reference triangle is the fixed one, then 
take the transform 
      of the curve $C$ under this map.
      
 The theorem is deduced from the fact that one can  
 arrive, after a finite number of iterations, at a 
 curve having only  
 ordinary multiple points as singularities. The fact 
that one 
 cannot obtain only ordinary double points as 
singularities 
 comes from step 2). 

 Once the elementary operation of blowing-ups of 
points was isolated, 
 it was possible to deduce immediately the 
following theorem from the 
 previous one:

\begin{theorem} \label{curvesurf}
    Let $C$ be a  reduced curve embedded in a 
smooth surface $S$. Then 
    $C$ can be resolved by a finite sequence 
of blowing-ups of points. 
    At each step of the algorithm, one simply 
blows-up the singular points of 
    the strict transform of $C$.  
\end{theorem}

Why is Theorem \ref{curvesurf} a consequence 
of Theorem \ref{quadr}? To
understand this,  
look at step 3) of the previous iteration. 
In it, the singular point
$c$ of the  
curve is blown up. At the same time other 
things happen to the plane, 
well described in Zariski's quotation at the 
end of subsection \ref{blow}. 
But if one concentrates one's attention in 
the neighborhood of the 
point $c$, the effect on the germ $(C,c)$ is the 
same as if one had 
only blown up that point. As it can be shown 
that any germ of reduced 
analytic curve embedded in a smooth surface 
can be embedded in the 
projective plane, one sees that the study 
done during the proof of 
Theorem \ref{quadr} allows one to prove 
Theorem \ref{curvesurf}. 

Let us be more explicit. The explanations 
which led to diagram
(\ref{ramnorm}) show that if $(C,c)\subset S$ 
is an irreducible germ
of curve in a smooth complex surface, then 
there are local coordinates
$(x,y)$ on $S$ centered at $c$, such that 
$(C,c)$ is given by a
parametrization of the form:
\begin{equation} \label{np}
  \left\{ \begin{array}{l}
               x=t^m\\
               y=\sum_{k\geq n}a_k t^k 
           \end{array} \right. .
\end{equation}
where $a_k \in \mathbb{C}, \: \forall \:k\geq n, 
a_n\neq 0$ and
$\mathrm{min}(m,n)$ is equal to the multiplicity 
of $C$ at $c$. A
parametrization of the form (\ref{np}) is 
called a \emph{Puiseux
  parametrization} or a \emph{Newton-Puiseux
  parametrization}. Such parametrizations are 
of the utmost
importance in the detailed study of singularities 
of plane curves (see
e.g. Brieskorn \& Kn{\"o}rrer \cite{BK 86}, 
Teissier
\cite{T 95}, Wall \cite{W 04}). 

It is possible to show that if $c$ is a singular 
point of $C$, that
is, if $m \geq 2$, then the choice of local 
coordinates can be done
such that $n >m$ and $n$ is not divisible by $m$. 
Then, as a
consequence of Proposition \ref{blowpoint}, one 
can show that the
strict transform of $C$ by the blowing-up of $c$ 
on $S$ can be
parametrized in suitable local coordinates by:
$$ \left\{ \begin{array}{l}
               x_1=t^m\\
               y_1=\sum_{k\geq n}a_k t^{k-m} 
           \end{array} \right. .$$
Continuing like this, we see that after 
exactly $[\frac{n}{m}]$
blowing-ups on the strict transform of $C$, 
one arrives at a strict
transform with multiplicity strictly less than $m$. 
Therefore,
\emph{multiplicity can be dropped by doing 
blowing-ups}. In the same
way, one can show that \emph{the intersection number 
of two germs of curves
embedded in a surface diminishes strictly after 
one blowing-up of
their intersection point}. The theorem is 
a direct consequence of these
two facts. 

\medskip
In the previous theorem it is not essential 
to suppose that 
the curve can be embedded, even locally, in 
a smooth surface. One has 
in general:

\begin{theorem}
   Let $(C,c)$ be a germ of reduced curve. Then 
it can be resolved by a finite 
   sequence of blowing-ups of points. 
At each step of the algorithm, one 
   simply blows-up the singular points of 
its strict transform.
 \end{theorem}

Let us sketch a more intrinsic proof (that is, 
which does not work
with local coordinates) for the case when the germ is
irreducible. Consider  the normalization morphism 
$(\overline{C}, \overline{c}) \stackrel{\nu}
{\rightarrow} (C,c)$. 
One has the inclusion of the corresponding 
local rings: 
$\mathcal{O}_{C,c} \subset \mathcal{O}_
{\overline{C}, \overline{c}}$. 
Denote by $F$ their common field of fractions. 
Denote by 
$(C_k, c_k)\stackrel{\pi_k}{\longrightarrow} 
(C,c), \: k \geq 1$ 
the composition of the first $k$ blowing-ups 
of the germ $(C,c)$ or 
of its strict transforms. Denote by  
$\mathcal{O}_k$ the local ring of the germ 
$(C_k, c_k)$.  
By Theorem \ref{exnorm}, the normalization 
morphism $\nu$ 
can be factored through the morphism $\pi_k$, 
which shows that 
one has the sequence of inclusions:
$$\mathcal{O}_{C,c}\subset \mathcal{O}_1 \subset
      \mathcal{O}_2 \subset \cdots \subset 
       \mathcal{O}_{\overline{C}, \overline{c}}\subset F.$$
 As $\nu$ is a finite morphism, one deduces that 
 $\mathcal{O}_{\overline{C}, \overline{c}}/ 
       \mathcal{O}_{C,c}$ is a finite dimensional 
  $\mathbb{C}$-vector space, which shows that 
one has to arrive at an 
  index $p \geq 1$ such that 
$\mathcal{O}_p = \mathcal{O}_{p+1}$. 
By Proposition \ref{blowpoint}, point 3), 
one sees
that $\mathcal{O}_{p+1}=
   \mathcal{O}_p[\frac{y_2}{y_1},...,
     \frac{y_r}{y_1}]$, where 
  $y_1,...,y_r$ are generators of the maximal ideal 
of the 
  local ring $\mathcal{O}_p$, chosen such that $y_1$ 
has 
  the smallest multiplicity when we look at 
the generators 
  as functions on the germ $(\overline{C}, 
\overline{c})$. 
  The equality $\mathcal{O}_p = \mathcal{O}_{p+1}$ 
  implies then that 
     $$\frac{y_2}{y_1},...,\frac{y_r}{y_1}\in 
\mathcal{O}_p$$
  which shows that the maximal ideal $(y_1, ..., y_r)
\mathcal{O}_p$ is principal, 
  and generated by $y_1$. But this shows that the local 
  ring $\mathcal{O}_p$ is regular, that is, the germ 
  $(C_p, c_p)$ is smooth. Again by Theorem 
  \ref{exnorm}, we deduce that $\pi_p$ is the 
normalization 
  morphism of $(C,c)$. As a consequence, one has 
  desingularized the germ $(C,c)$ after $p$ 
blowing-ups.
  
  When $(C,c)$ is not irreducible, the 
  total ring of fractions of $\mathcal{O}_{C,c}$ 
is no more 
  a field, but a direct product of fields. 
At some steps 
  of the blowing-ups irreducible components may be 
  separated, but the overall analysis remains 
the same. 

For a careful proof written in the language 
of commutative algebra 
and for many details on abstract singularities 
of not necessarily plane 
curves, one can consult Castellanos \& 
Campillo's book  \cite{CC 05}. 
\medskip

The next theorem is a generalization of Theorem 
\ref{genusmooth}. It shows how the finite 
dimension of the 
quotient $\overline{\mathcal{O}}_{C,c}/
\mathcal{O}_{C,c}$ used in the 
previous proof appears in the computation of 
the genus of 
the normalization. By $\overline{\mathcal{O}}_{C,c}$ 
we denote the 
integral closure of $\mathcal{O}_{C,c}$ in its total 
ring of 
fractions, that is (see Theorem \ref{transcalg}), 
the direct sum of
the local rings of the normalization  
$\overline{C}$ at the preimages of the point $c$. 

\begin{theorem} \label{genusing}
    Let $C$ be a reduced algebraic curve inside 
$\C\mathbb{P}^2$. 
    Then the genus of its normalization $\overline{C}$ 
is equal to 
    $\frac{(d-1)(d-2)}{2}- \sum\delta(C,c)$, 
where the sum is done 
    over the singular points of $C$, and at 
such points 
    $\delta(C,c) := \mathrm{dim}_{\C} \: 
    \overline{\mathcal{O}}_{C,c}/\mathcal{O}_{C,c}$. 
 \end{theorem}

For many more details about possibly singular 
plane curves, we 
recommend the leisurely introduction done 
in Brieskorn \& 
Kn{\"o}rrer \cite{BK 86}.

\medskip
\subsection{Embedded resolution of plane curves}$\:$
\medskip

Let us consider again Theorem \ref{curvesurf}.  
If $\pi$ denotes the 
composition of the blowing-ups which resolves 
the curve $C$, 
one knows by the definition of resolution that 
the strict transform 
$C'_{\pi}$ is smooth. But the total transform 
$\pi^{-1}(C)$ has not 
necessarily only normal crossings. Nevertheless, 
by blowing-up 
more, one can arrive at an embedded resolution:

\begin{theorem} \label{algemb}
   Let $C\hookrightarrow S$ be a reduced curve 
embedded in a smooth 
complex analytic surface. Start from the 
identity morphism 
$S_0=S\stackrel{\pi_0}{\rightarrow}S$. Then the
following algorithm stops after a finite number 
of steps:

$\bullet$ If $S_k \stackrel{\pi_k}{\rightarrow}S$ 
is given and the
total transform $\pi_k^{-1}(C)$ has more 
complicated singularities
than  normal crossings inside
the surface $S_k$, then blow-up each point of 
$\pi_k^{-1}(C)$ at which 
its irreducible components do not cross normally. 
The composition of
$\pi_k$ and of these blowing-ups is 
by definition $S_{k+1}
\stackrel{\pi_{k+1}}{\rightarrow}S$. 

$\bullet$ If $\pi_k^{-1}(C)$ has normal crossings inside
the surface $S_k$, then STOP. 
\end{theorem}

Moreover, the embedded resolution obtained in 
this way can be
distinguished among all embedded resolutions 
by a \emph{minimality
  property}, to be compared with the one stated 
in Theorem \ref{minresurf}:

\begin{theorem}
   Let $S_{min} \stackrel{\pi_{min}}{\longrightarrow}S$ 
be the
   embedded resolution of the curve $C$ obtained 
by the algorithm 
   \ref{algemb}. Then $\pi_{min}$ is minimal among 
all the embedded
   resolutions of $C$, in the sense that any 
other resolution 
   $S' \stackrel{\pi'}{\rightarrow}S$ factorizes 
as $S' \stackrel{\psi}{\rightarrow}S_{min} 
  \stackrel{\pi_{min}}{\longrightarrow}S$, where 
$\psi$ is a
   composition of blowing-ups of points. 
\end{theorem} 

One can use the previous theorem as a way 
to analyze the structure of 
a singular point of a curve embedded in a smooth 
surface. More 
precisely, one can look at various aspects of 
the resolution 
$\pi_{min}$ and of the sequence of blowing-ups 
leading to it:
\makeatletter
\renewcommand{\theenumi}{\alph{enumi}}
\renewcommand{\labelenumi}{(\theenumi)}
\makeatother
\begin{enumerate}
    \item \label{multseq}
        The sequence of multiplicities of 
the strict transforms of the germ.
    \item \label{dualgraphmult}
         The dual graph of the total transform 
of $(C,c)$, each vertex 
         being decorated by the order 
of vanishing of the preimage 
         of the maximal ideal $\mathcal{O}_{S,c}$ 
on the corresponding 
         component.
    \item \label{dualgraphint}
        The dual graph of the total transform of 
$(C,c)$, each vertex 
         being decorated by the self-intersection 
number of the 
         corresponding component. 
    \item \label{proximity}
          A graph which represents the strict 
transforms of the 
          components of the curve $(C,c)$ 
at each step of the 
          process of blowing-ups, and whose edges 
          are drawn in such a way as to remember 
if the strict 
          transform passes or not through a 
smooth point of the 
          exceptional locus.
 \end{enumerate} 
 
It can be shown that all these encodings  
are equivalent. They are 
also equivalent with information readable 
on Newton-Puiseux
parametrizations of the germ  $(C,c)$ (that is, 
parametrizations of
the type (\ref{np})).  Furthermore, all 
the encodings describe 
completely the \emph{embedded topological type} 
of the germ. 
For the state of the art around 2004 on 
the relations with the
embedded topology of germs, see Wall \cite{W 04}. 

The comparison between information readable on 
Newton-Puiseux 
expansions and on sequences of blowing-ups 
(in fact quadratic 
transformations, as we explained before) 
seems to have been 
started by Max Noether in \cite{N 90}. 
It was carefully explored 
by Enriques \& Chisini \cite{EC 18}, 
who introduced the 
viewpoint (\ref{proximity}) of the previous list. 
A recent textbook
emphasizing  the usefulness of 
such graphs (called nowadays \emph{Enriques diagrams}) 
in the study of plane curve singularities is 
Casas-Alvero \cite{CA 00}. 

A detailed comparison between the viewpoint 
(\ref{dualgraphmult}) 
and information readable on Newton-Puiseux 
expansions was 
done in Garc{\'\i}a Barroso \cite{GB 96}.

The four viewpoints (\ref{multseq})-(\ref{proximity}) 
are also 
compared in Campillo \& Castellanos \cite{CC 05}, 
where some 
of them are extended to arbitrary reduced germs, 
not 
necessarily embeddable in smooth surfaces.

A common aspect of all the comparisons is the use 
of expansions of
rational numbers into continued 
fractions. In \cite{PP 06} we made a detailed 
study 
of the convex geometry lying behind the use 
of continued fractions, and 
of its applications via toric geometry to the study 
of singularities
of curves and surfaces.

\medskip
\section{Resolution of surface singularities by Jung's  method} \label{Jung}

\medskip
\subsection{Strategy}\label{strat} $\:$
\medskip

Our aim in this section is to prove the following:

\begin{theorem} \label{resurf}
   Any reduced complex surface admits a resolution.
\end{theorem}

Unlike for the case of plane curves, 
the normalization morphism is no
longer always a resolution. 
Indeed, the world of \emph{normal} surface
singularities is huge and fascinating. 
For example, all isolated surface
singularities of complete intersections are normal. 
But as the
normalization morphism exists and it is 
an isomorphism over the smooth
locus, one can reduce the proof of 
Theorem \ref{resurf} to the task of 
proving it for normal surfaces, 
although this is not much of a
simplification. 

Many methods to prove Theorem \ref{resurf} 
have been proposed,
sometimes for 
special types of surfaces (projective surfaces in
$\mathbb{P}^3$, arbitrary projective surfaces, 
algebraic surfaces,
analytic surfaces, arithmetic surfaces, etc.) 
For the first approaches on
the problem, one can consult Gario's papers 
\cite{G 89}, \cite{G 91}. 
For the state of the art around 1935, 
see Zariski's book \cite{Z 35}. For the
progresses made up to the year 2000, 
one can consult Hauser's index
\cite{H 00}. We 
recommend also Lipman's survey \cite{L 75},  
Cutkosky's book \cite{C  04} and Koll\'ar's book 
\cite{K 07}. In these last two references, one 
may find detailed proofs of the existence 
of resolution of singularities of complex 
algebraic varieties in arbitrary dimension.

We will now explain one of the methods, 
 usually known  as
\emph{Jung's method}. It is probably the most 
amenable one to
computations by hand on 
examples defined by explicit equations. 
For special types of
singularities, other methods could be more suitable. 
For example, if
the surface admits a $\mathbb{C}^*$-action, 
then one knows how to
describe equivariantly a resolution 
(see Orlik \& Wagreich \cite{OW
  71} or M{\"u}ller \cite{M 00}). 

From a very general viewpoint, one could express 
the fundamental idea of this method by the injunction:
\begin{equation} \label{injung}
     \begin{array}{c} 
         \mbox{\emph{In order to represent 
an object as the image of
         a simpler one}} \\  
          \mbox{\emph{first choose an image 
of the object}} \\
          \mbox{\emph{then simplify this image.}}
     \end{array} 
\end{equation}

Given a reduced complex analytic surface, 
Jung's method consists in
ana\-lysing its structure in a neighborhood  
of one  of its singular points by projecting it 
to a plane and by
considering an embedded  
resolution of the discriminant curve.  
It was introduced by Jung
\cite{J 08} as a way to  
uniformize locally a surface, and extended by 
Walker \cite{W 35} in
order to prove  
resolution of singularities of algebraic surfaces. 
This second paper
was considered by  
Zariski \cite[Chapter 1]{Z 35} to be the first 
complete proof of the
resolution of singularities  
for surfaces. Hirzebruch \cite{H 53} used again 
the method in order to
prove the  
resolution of singularities for complex analytic 
surfaces (for an
excellent summary of  
Hirzebruch's work on singularities we refer 
to Brieskorn \cite{B 00}). Here we
will explain Hirzebruch's  
proof (see also Laufer \cite[Chapter  II]{L 71}
and Lipman \cite{L 75}).

Before starting our explanation, we would like 
to emphasize that
Hirzebruch's motivation was to extend to the 
case of $2$ complex
variables Riemann's construction of a smooth (real) 
surface associated
to a multivalued analytic function of one variable. 
We quote from the
introduction of \cite{H 53}:

\begin{quote}
    {\small The <<algebroid>> function elements 
of a multivalued
      function $f(z_1, z_2)$ defined in a complex 
manifold (with two
      complex dimensions) can be easily associated 
to the points of a
      Hausdorff space of dimension 4, which covers 
part of $M$, and
      which we will call \emph{the Riemann domain} 
of the function
      $f$. But this Riemann domain is not in general 
a topological
      manifold....} 
\end{quote}

The main steps of Jung's method 
of resolution of reduced
complex analytic surfaces  are:
\makeatletter
\renewcommand{\theenumi}{\Alph{enumi}}
\renewcommand{\labelenumi}{(\theenumi)}
\makeatother
\begin{enumerate}
   \item Take a germ of the given surface 
and consider a finite
       morphism to a germ of   smooth surface.
       
   \item Consider an embedded resolution of 
the discriminant curve of
             this morphism,  
             and pull-back the initial germ by 
this resolution morphism.
             
   \item Normalize the surface obtained by 
this pull-back.
   
   \item Resolve explicitly the singularities of 
the new normal
             surface, by using the fact  
             that they admit a finite morphism 
to a smooth surface
             whose discriminant curve  
             has normal crossings.
             
   \item Glue together all the previous constructions 
to get a global
   resolution of the initial surface.     
\end{enumerate}

\medskip

Let us explain the previous steps with more details. 
We emphasize that the surface 
$S$ is not supposed to be normal. 
\medskip

(A) Let $S$ be the given reduced surface. Consider 
one of its points $s \in S$,
and the germ $(S,s)$.  Denote  
by $(S,s)\stackrel{\alpha}{\rightarrow}  (R,r)$ 
a finite morphism,
where $(R,r)$ is a germ of smooth surface.  
Consider its discriminant locus $\Delta(\alpha) 
\subset R$. Three
different situations may occur: either it is empty, 
either it is the point
$r$, either it is a germ of curve.  
\medskip

(B) When $\Delta(\alpha)$ is empty, $\alpha$ 
is a local isomorphism,
which shows that $S$ is smooth at $s$. 

When $\Delta(\alpha)=r$, the normalization 
morphism $\coprod
(\overline{S}_i, \overline{s}_i) 
\stackrel{\chi}{\rightarrow} (S,s)$,
whose total space is a multigerm, is also a 
resolution. Indeed, each
restriction $(\overline{S}_i, \overline{s}_i)
\stackrel{\alpha \circ \chi}{\longrightarrow} (R,r)$ 
is unramified
outside $r$, and as $\pi_1(V\setminus r)=0$ for 
each polydisc
representative of $R$, one sees that a finite 
representative of
$\alpha\circ \chi$ is a trivial covering 
over $V \setminus r$. The
generalized Riemann extension theorem \ref{rext2} 
implies that $\alpha
\circ \chi$ is an isomorphism, which shows that 
each $(\overline{S}_i, \overline{s}_i)$ is
smooth.

Let us suppose finally that 
$(\Delta(\alpha),r)\subset 
(R,r)$ is a germ of curve. Take an embedded 
resolution $(\tilde{R}, E)
\stackrel{\psi}{\rightarrow} (R,r)$ of the germ, 
then 
the pull-back $\psi^*(\alpha)$ of the morphism 
$\alpha$. One gets like
this the following diagram of analytic morphisms:
\begin{equation} \label{pullback}
\xymatrix{
    (\tilde{S}, F)
         \ar[d]_{\psi^*(\alpha)} 
             \ar[rr]^-{\alpha^*(\psi)} & 
       &  (S, s)
                 \ar[d]^{\alpha} \\
    (\tilde{R}, E) \ar[rr]^{\psi}& & (R,r)  }
\end{equation}

  As the mapping $\psi$ is a proper modification, 
one deduces that
  $\alpha^*(\psi)$ is also a proper modification.

By construction, $\psi^{-1}(\Delta(\alpha))$ is 
a curve with normal crossings
inside $\tilde{R}$ and $\psi^{-1}(\Delta(\alpha))=E \cup
\Delta(\alpha)_{\psi}'$, where $\Delta(\alpha)'_{\psi}$ 
denotes the strict
transform of $\Delta(\alpha)$ by the 
modification $\psi$. 

The discriminant locus $\Delta(\psi^{*}(\alpha))$ of 
the morphism
$\psi^*(\alpha)$ is contained in 
$\psi^{-1}(\Delta(\alpha))$. Therefore
the germs of the surface $\tilde{S}$ at the points 
of $F=
(\alpha^*(\psi))^{-1}(s)$ have a special property: 
they can be
projected by finite morphisms (the localization of
$\psi^*(\alpha)$) to a smooth surface, such that 
the discriminant
locus has normal 
crossings. In many ways such germs are much more 
tractable than
arbitrary germs of surfaces, that is why they 
received a special name:

\begin{definition} \label{qobj}
Let  $(\mathcal{X},x)$ be a germ of reduced 
equidimensional complex
surface.  The germ $(\mathcal{X},x)$ is called 
\textbf{quasi-ordinary} if there exists a finite 
morphism $\phi$ from 
$(\mathcal{X},x)$ to a germ of smooth surface, 
whose discriminant locus is contained in a curve with 
normal crossings. Such a morphism $\phi$ is also 
 called 
\textbf{quasi-ordinary}.
\end{definition}

An example of quasi-ordinary germ is the germ at 
the origin of 
 Whitney's umbrella defined by equation 
(\ref{Whitneyeq}). 
A quasi-ordinary morphism associated to it is 
the restriction of the 
canonical projection of the ambient space to  
the plane of 
coordinates $(y,z)$. Please contemplate how this 
is visible 
in Figure 4, where the discriminant curve is 
the union of two 
lines, one being the projection of the singular 
locus and the other 
one being the apparent contour. At this point, 
recall also the 
comments about the drawing of surfaces from 
the beginning 
of Subsection \ref{finmor}.

The name ``quasi-ordinary'' was probably introduced 
in reference to the
previously named ``ordinary'' singularities. 
In \cite[page 18]{Z 35},
Zariski says that a 
surface in the projective space has an ordinary 
singularity at a point
if either it is locally isomorphic at this point 
to a singular normal crossings
divisor, or to an ``ordinary cuspidal point'', 
defined geometrically. In
fact these last germs are isomorphic to the germ 
at the origin of 
Whitney's umbrella. 

 In many respects, quasi-ordinary germs are more 
amenable to  study 
 than arbitrary 
 singularities, because one can extend to them 
by analogy many
 constructions done first for curves. For example, 
in what concerns
 resolution of singularities, Gonz{\'a}lez P{\'e}rez 
\cite{GP 03} gave two 
 methods for finding embedded resolutions of 
quasi-ordinary germs of
 hypersufaces in arbitrary dimensions, by developing 
a method  analogous
 to the one proposed for the case of curves 
by Goldin \&
 Teissier \cite{GT 00}. As 
 an introduction to quasi-ordinary singularities 
in arbitrary
 dimensions, we recommend Lipman's foundational 
work \cite{L 88}.

\medskip

(C) Coming back to the diagram (\ref{pullback}), 
let us normalize the
    surface $\tilde{S}$. Denote by 
$\overline{\tilde{S}}
    \stackrel{\nu}{\rightarrow} \tilde{S}$ 
the normalization
    morphism. One gets the diagram:

\begin{equation} \label{normback}
 \xymatrix{(\overline{\tilde{S}}, G)
\ar[ddr]_{\overline{\alpha}} 
             \ar[dr]_{\nu} \ar[drrr]^{\overline{\psi}} 
   & & & \\
    & (\tilde{S}, F)
         \ar[d]^{\psi^*(\alpha)} 
             \ar[rr]_{\alpha^*(\psi)} & 
       &  (S, s)
                 \ar[d]^{\alpha} \\
    & (\tilde{R}, E) \ar[rr]^{\psi}& & (R,r)  }
\end{equation}

By definition, $\overline{\alpha}:= \psi^*(\alpha) 
\circ \nu$ and 
   $\overline{\psi}:= \alpha^*(\psi) \circ \nu$. 
As  normalization
   morphisms are finite modifications, we see 
that $\overline{\alpha}$ is
   finite and $\overline{\psi}$ is a modification. 
Moreover, the
   discriminant locus $\Delta(\overline{\alpha})$ of
   $\overline{\alpha}$ is contained in the 
discriminant locus of
   $\psi^*(\alpha)$, which shows that 
$\Delta(\overline{\alpha})$ has
   again normal crossings. Therefore, the 
singularities of
   $\overline{\tilde{S}}$ are still more special 
than those of
   $\tilde{S}$. They are the so-called 
   \emph{Hirzebruch-Jung singularities}:

\begin{definition} \label{defhj}
  A \textbf{Hirzebruch-Jung} germ (or singularity) 
of complex surface is a
  \emph{normal} quasi-ordinary germ of surface.  
\end{definition}

\medskip
(D) We explain in the next subsection how one can use 
    Definition \ref{defhj} directly  in order 
to give an explicit resolution of
    any Hirzebruch-Jung singularity. For 
the moment, please accept the 
    fact that Hirzebruch-Jung singularities 
admit resolutions. 

    As the surface
    $\overline{\tilde{S}}$ has only this special 
kind of singularities, one
    sees that it can be resolved. Denote by
    $T\stackrel{\rho}{\rightarrow}
\overline{\tilde{S}}$ a resolution
    of $\overline{\tilde{S}}$. 
Then $\alpha^*(\psi)\circ\nu\circ\rho$ is
    a modification of $(S,s)$, being a composition 
of three
    modifications. Moreover, its source is smooth 
and it is an
    isomorphism over the smooth locus of $S$, 
which shows that it is a
    resolution of the germ $(S,s)$. 

\medskip
(E) Once the germ $(S,s)$ was fixed, the finite 
morphism $\alpha$ used to 
project it on a smooth surface was arbitrary. 
One can choose then 
a representative $U \stackrel{\alpha}{\rightarrow} V$ 
of the germ of 
morphism $\alpha$ which is a finite morphism of 
analytic surfaces, and 
such that the germ $(\Delta(\alpha),r)$ admits 
a closed representative 
in $V$ which is smooth outside $r$. Then the 
composed morphism 
$\alpha^*(\psi)\circ\nu\circ\rho$ is by 
construction a normalization of 
$U \setminus s$ (recall that $s$ is not 
necessarily an isolated
singular point of $S$). Indeed, its source 
is smooth, it is proper and
bimeromorphic as a composition of 
proper and bimeromorphic morphisms, and its 
only possible fiber of positive 
dimension lies over the point $r$. We use 
then Theorem \ref{exnorm} to
complete the argument.

This shows that the set of points of the 
surface $(S,s)$ which do not have 
neighborhoods resolved by the normalization of $S$ 
is \emph{discrete}.
If $\Sigma(S)$ is this set, choose neighborhoods 
$(U_s)_{s \in \Sigma(S)}$ which are pairwise disjoint 
and which are at
the same  
time sources of finite morphisms as explained in 
the previous
paragraph. Apply then  
the previous process for each one of them, getting 
like this 
resolution morphisms $T_s \rightarrow U_s$ for 
any 
$s \in \Sigma(S)$. If one considers also the 
normalization morphism 
of $S \setminus \Sigma(S)$, all these modifications 
of open sets which
form 
a covering of the surface $S$ 
\emph{agree on overlaps}. This implies 
that they can be \emph{glued} into a resolution 
of the entire surface  $S$. 

\medskip
Examples of applications of this method are given 
in 
Laufer \cite{L 71}, L{\^e} \& Weber \cite{LW 00} 
and N{\'e}methi 
\cite{N 99}. In this last reference are described 
the dual 
graphs of the resolutions obtained by this method 
for the 
germs of surfaces defined by equations of the form
$z^n + f(x,y)=0$, where $f$ is reduced. 
The description is 
done in terms of the embedded resolution of the curve 
$f(x,y)=0$ and the integer $n\geq 1$.

We see that the application of the injunction 
(\ref{injung}) depends
heavily on the possibility to do an embedded 
resolution of curves
contained in smooth surfaces. Nevertheless, 
this  does not lead to an
embedded resolution of a surface, only to 
an abstract one. If the
method  could
be adapted to give an embedded resolution, 
then the same strategy 
would prove the resolution of 3-folds, and 
then one could try to
get on these lines an inductive resolution 
in all dimensions. Untill
now, this strategy was not succesfull 
(see Problem (\ref{first}) in Section \ref{open}).

\medskip
\subsection{Resolution of Hirzebruch-Jung 
singularities}$\:$
\medskip

The aim of this section is to explain 
a proof of the following theorem:

\begin{theorem} \label{resolhj}
    Let $(S,s)$ be a Hirzebruch-Jung singularity. 
Then it can be resolved.
    Moreover, the exceptional divisor of its 
minimal resolution has normal 
    crossings, its components are smooth rational 
curves and its 
    dual graph is a segment.
\end{theorem}

The method we will present is a generalization of 
the one used to construct
topologically  
a normalization morphism for an irreducible 
germ of curve 
(see the explanations which 
precede diagram (\ref{ramnorm})). 
We saw there that one could work in
convenient local coordinates with a morphism defined by
$x_1=t^m$. Similarly, one can resolve Hirzebruch-Jung 
surface
singularities using only morphisms defined by 
{\em monomials} in two
variables. There is a branch of algebraic 
geometry which studies
intrinsically such morphisms, called 
\emph{toric geometry}. In what
follows we will explain how it appears in our context. 
For an
introduction to toric geometry, one may consult 
Fulton \cite{F 93}.
\medskip

Let $(S,s)$ be a Hirzebruch-Jung singularity and 
$(S,s)\stackrel{\alpha}{\rightarrow}  (R,r)$ be 
a finite morphism whose discriminant locus has normal 
crossings. 
Choose local coordinates in the neighborhood of 
the point $r \in R$ such that 
the discriminant curve is contained in the union 
of the 
coordinate axis. 
Therefore, from now on we suppose that 
$R=\mathbb{C}^2$ and that 
$\alpha$ is unramified over $(\mathbb{C}^*)^2$. 

Choose a finite representative 
$U\stackrel{\alpha}{\rightarrow}  V$ where 
$V$ is a polycylinder in the coordinates $(x,y)$. 
Denote by $V^*$ the 
complement of the axis of coordinates and by $U^*$ 
the preimage  
$\alpha^{-1}(V^*)$. Then the restricted morphism
$U^*\stackrel{\alpha}{\rightarrow}  V^*$  
is a finite (unramified) covering with connected 
source (because 
$(S,s)$ was supposed to be normal). Consider the 
associated morphism of fundamental groups
$\pi_1(U^*)  \stackrel{\alpha_*}{\rightarrow}  
\pi_1(V^*)$. As $V$ was 
chosen to be a polycylinder, $V^*$ is the product 
of two pointed 
discs, which shows that  $N:=\pi_1(V^*)$ is 
 a free abelian group of
rank 2.  Therefore  
all its finite index subgroups are also 
free abelian of rank 2.  In
particular this holds for   
$N(\alpha):= \alpha_*(\pi_1(U^*))$. Notice that 
the abelianity of
those groups implies that we 
do not have to worry about base points. 

Let us look  at the multiplicative semigroup of 
Laurent monomials
$x^ay^b$ which are regular in a  
neighborhood of the origin of $\mathbb{C}^2$, 
which means that $a \geq
0, b\geq 0$. One can think about them as 
functions $(\C^*)^2
\rightarrow \C^*$ and  
restrict them to  the loops which represent 
the elements of 
$N(\alpha)$. One associates like this to 
each pair $(a,b)$ a map from
an oriented  
circle to $\mathbb{C}^*$, map whose degree (its linking 
coefficient with the origin) is a well-defined 
integer. Consider the
elements of $N(\alpha)$ whose associated  
degree is \emph{non-negative}. They form a 
sub-semigroup 
of $(N(\alpha),+)$.  Denote  by $\sigma$ the 
closed convex 
cone generated by its elements inside the 
associated real vector space 
$N(\alpha)_{\R}:= N(\alpha)\otimes_{\Z}\R$. 
As $N(\alpha)$ has finite
index in $N$, one has a canonical 
identification $N(\alpha)_{\R} =
N_{\R}$. Therefore, $\sigma$ may also be seen as 
a cone in
$N_{\R}$. It is precisely the first quadrant!

The pair 
$(N(\alpha), \sigma)$ determines a two-dimensional 
normal affine toric surface $\mathcal{X}( N(\alpha), 
\sigma)$ and 
the inclusion $N(\alpha) \hookrightarrow N$ induces  
 a canonical
toric morphism:  
$$\mathcal{X}( N(\alpha), \sigma) \stackrel{\gamma_
     {N:N(\alpha)}}{\longrightarrow} \mathbb{C}^2
= \mathcal{X}( N, \sigma) .$$
One sees then that the morphism of fundamental 
groups 
induced over $(\C^*)^2$ has $N(\alpha)$ as its image. 
From this, one concludes 
that $\gamma_{N:N(\alpha)}$ can be lifted over 
$(\C^*)^2$ to a
morphism with target space  
$(S,s)$. By construction, this lift  is bounded in 
a neighborhood of the 
special point $\underline{0}$ of $\mathcal{X}( 
N(\alpha), \sigma) $ (its unique
closed orbit under the torus action). As this 
last variety 
is normal, one deduces that it can be extended 
to a morphism 
$\mu$ defined on a neighborhood of 
$\underline{0}\in \mathcal{X}( N(\alpha),
\sigma)$. Thus, one gets a  diagram:
\begin{equation} \label{isohj}
         \xymatrix{
           (\mathcal{X}( N(\alpha), \sigma) ,
\underline{0}) \ar[r]^-{\mu} 
             \ar[dr]_{\gamma_{N:N(\alpha)}}& (S,s) 
\ar[d]^{\alpha}\\
             &  (\mathbb{C}^2,0) }
  \end{equation}
By construction,  $\mu$ is a finite morphism which is 
an isomorphism 
over $(\mathbb{C}^*)^2$. As both its source and  
target are normal, one 
deduces that $\mu$ is an isomorphism.  This shows 
that a Hirzebruch-Jung 
singularity is analytically isomorphic to a germ 
of toric
surface.            
    
  But the singularities of normal toric surfaces 
admit explicit 
  minimal resolutions, which can be deduced from 
the geometry of the pair 
  $(N(\alpha), \sigma)$. They verify the properties 
listed in 
  Theorem \ref{resolhj}. 
  \medskip

  When Hirzebruch did his work \cite{H 53}, 
toric geometry did not exist. 
  Nevertheless, he gave an explicit  resolution of 
Hirzebruch-Jung 
  singularities in a way which nowadays can be 
recognized to be toric,
  by gluing affine planes through monomial maps. 
We recommend 
  Brieskorn's article \cite{B 00} for comments 
on this approach, as well as on 
  other contributions by Hirzebruch to singularity 
theory. 
  
  It was one of our contributions to the study of 
Hirzebruch-Jung
  singularities to  
  construct  the affine toric surface 
$\mathcal{X}( N(\alpha), \sigma)$. 
  Our motivation was to be able to compute 
the normalizations of 
  explicit quasi-ordinary singularities. 
We showed that in arbitrary 
  dimensions the normal quasi-ordinary 
singularities could be 
  characterized as the germs of normal affine 
toric varieties defined 
  by a simplicial cone.  Moreover, we gave an 
algorithm of 
  normalization for hypersurface quasi-ordinary 
germs in arbitrary 
  dimensions (see \cite{PP 05}). 

As a particular case of this normalization algorithm, 
one gets the
following lemma, which is needed when one applies 
steps (C) and (D) described in Section \ref{strat} 
to concrete
examples (as the ones presented in Laufer \cite{L 71}, 
L{\^e} \& Weber
\cite{LW 00} and N{\'e}methi  \cite{N 99}): 

\makeatletter
\renewcommand{\theenumi}{\arabic{enumi}}
\renewcommand{\labelenumi}{(\theenumi)}
\makeatother

\begin{lemma}
  Let $(S,s)\hookrightarrow \mathbb{C}^3_{x_1, x_2, y}$ 
be the
  quasi-ordinary irreducible singularity corresponding 
to the
  algebraic function with two variables $y:=
  x_1^{\frac{p_1}{q_1}}x_2^{\frac{p_2}{q_2}}$, where
  $\mathrm{gcd}(p_1, q_1)= \mathrm{gcd}(p_2, q_2)=1$ 
and $p_1, q_1,
  p_2, q_2\in \mathbb{N}^*$. Denote 
  $d:= \mathrm{gcd}(q_1, q_2), \: j_1=\frac{q_1}{d}$ 
and let $k_1\in
  \{0,1,...,q_1-1\}$ be the unique number in this 
set which satisfies
  the congruence equation $k_1p_1 +j_1p_2\equiv 0 \: 
(\mathrm{mod}\:
  q_1)$. Denote also:
  $$q_1':= \frac{q_1}{\mathrm{gcd}( q_1, k_1)}, 
\: k_1':=
  \frac{k_1}{\mathrm{gcd}( q_1, k_1)}.$$
  Consider $\frac{q_1'}{k_1'}=b_1-1/
(b_2-1/(\cdots-b_r))$, the
  decomposition of $\frac{q_1'}{k_1'}$ as 
a Hirzebruch-Jung continued
  fraction (that is,  
 $b_i \geq 2, \: \forall\: i\in \{ 1,...,r\}$).
  Then:
\begin{enumerate}
   \item $(S,s)$ has the same normalization 
over $\mathbb{C}^2_{x_1,
    x_2}$ as the surface $(S',s)$ defined by 
the algebraic function 
   $$ x_1^{\frac{1}{q_1'}}x_2^{\frac{q_1'-k_1'}{q_1'}}.$$
 
   \item The
 total transform  of the function $(x_1x_2)|_{(S,s)}$ 
 by the minimal resolution of $(S,s)$ has 
a dual graph as
 drawn in Figure 7, where $\{x_i=0\}'$ denotes 
the strict transform of
 $x_i=0$. 
\end{enumerate}
\end{lemma}
 
The important point to notice is the way 
the strict transforms of the
germs of curves defined by $x_1=0$ and $x_2=0$ 
intersect the
exceptional divisor of the minimal resolution 
of $(S,s)$. Be careful
not to permute $x_1$ and $x_2$!

{\tt    \setlength{\unitlength}{0.92pt}}
 \begin{figure} 
 \epsfig{file=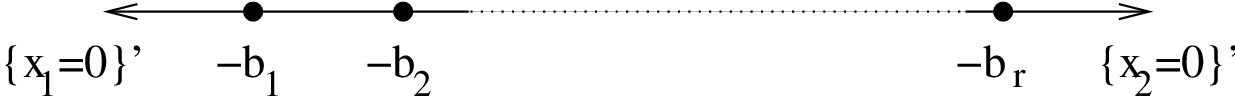, height= 10mm}
  \caption{Total transform of $x_1x_2=0$}
 \end{figure}
\medskip

We conclude this section with a theorem 
which characterizes
Hirzebruch-Jung surface 
singularities from many different viewpoints.  
References can be found in
\cite{BHPV 04}, \cite{PP 05} and \cite{PP 06}.

\begin{theorem} \label{caracthj}
   Hirzebruch-Jung surface singularities can 
be characterized
   \emph{among normal singularities} by 
the following 
   equivalent properties:
   \begin{enumerate}
       \item They are quasi-ordinary.
       \item They are singularities of toric surfaces.
       \item They are cyclic quotient singularities. 
       \item The exceptional divisor of 
their minimal resolution has 
                 normal crossings, its components 
are smooth rational 
                 curves and its dual graph is 
a segment.
       \item Their link is a lens space.
       \item The fundamental group of their link 
is abelian. 
   \end{enumerate}              
\end{theorem}

\medskip
\section{Open problems} \label{open}

\begin{enumerate}

  \item \label{first} 
Adapt Jung's method to get {\em embedded} 
resolution of germs
    of surfaces in $\C^3$. 

  \item Use Jung's method to get obstructions on 
the topology of
    germs of surfaces with isolated singularities 
in $\C^3$.  

  \item If $(S,s)$ is a germ of normal surface  
and $(S,s)
    \stackrel{\alpha}{\rightarrow} (\C^2,0)$ is 
a finite morphism, one
    gets by Jung's method an associated resolution 
of $(S,s)$. In
    general, one gets more components of 
the exceptional divisor than
    in the minimal resolution. Denote 
by $\mathrm{md}(S,s)$ the minimum number
    of supplementary components when 
one varies $\alpha$, the germ
    $(S,s)$ being fixed (`\emph{md}' being 
the initials of `\emph{minimal difference}'). 
Is $\mathrm{md}(S,s)$ bounded 
from above when one varies
    $(S,s)$ among the normal germs with fixed 
topology? 

    By  construction, $\mathrm{md}(S,s)$ attains 
a minimal value when one varies
    $(S,s)$ like this. Compute it in terms of 
the weighted dual graph
    of the minimal good resolution of $(S,s)$ 
(which encodes the
    topology of $(S,s)$, as ensured by a theorem 
of Neumann).

\end{enumerate}

\medskip

\textbf{Acknowledgements:} I am grateful to 
Angelica Cueto, Eleonore Faber, 
Camille Pl{\'e}nat, Jawad Snoussi and 
 Dmitry Stepanov for their pertinent
remarks on a previous version of this paper.


\medskip

\end{document}